\documentclass[11pt, reqno]{amsart}
\usepackage{graphicx,xcolor,amsmath,amssymb,comment,bm,hyperref}
\usepackage[shortlabels]{enumitem}

\newtheorem{thm}{Theorem}[section]

\newtheorem{lemma}[thm]{Lemma}

\newtheorem{coro}[thm]{Corollary}
\newtheorem{obs}[thm]{Observation}
\newtheorem{remark}[thm]{Remark}
\newtheorem{prop}[thm]{Proposition}

\newcommand{\bs}{\backslash}

\usepackage{tikz}
\usetikzlibrary{decorations.pathmorphing,decorations.pathreplacing,decorations.markings, decorations.shapes}
\tikzset{snake it/.style={decorate, decoration=snake}}
\tikzset{decorate sep/.style 2 args=
	{decorate,decoration={shape backgrounds,shape=circle,shape size=#1,shape sep=#2}}}

\newcommand{\dec}{\draw[decorate sep={.3mm}{1mm},fill]}

\begin{document}

	\title[Unavoidable Immersions]{Unavoidable immersions of $4-$ and $f(t)-$edge-connected graphs}

	\author{Guoli Ding}
	\address{Mathematics Department\\
		Louisiana State University\\
		Baton Rouge, Louisiana}
	\email{ding@math.lsu.edu}
	
	\author{Brittian Qualls}
	\address{Mathematics Department\\
		Louisiana State University\\
		Baton Rouge, Louisiana}
	\email{bquall3@lsu.edu}

	\date{\today}

	\begin{abstract}
		In this paper we prove that every sufficiently large 4-edge-connected graph contains the double cycle, $C_{2,r}$, as an immersion. In proving this, we develop a new tool we call a ring-decomposition. We also prove that linear edge-connectivity implies the presence of a $C_{t,r}$ immersion in a sufficiently large graph, where $C_{t,r}$ denotes the graph obtained from a cycle on $r$ vertices by adding $(t-1)$ edges in parallel to each existing edge; this result is an edge-analogue of a result of B\"{o}hme, Kawarabayashi, Maharry, and Mojar. We then use the latter result to provide an unavoidable minor theorem for highly connected line graphs.
	\end{abstract}

	\maketitle

	\section{Introduction}
	
	In this paper, all graphs are finite and may have multiple edges; since loops do not add to the theory of graph immersions, we will consider all graphs in this paper to be loopless. A graph $H$ is called a
	\textit{minor} of a graph $G$ if a graph isomorphic to $H$ can be obtained from a subgraph of $G$ by contracting edges. A graph $H$ is called a \textit{topological minor} of a graph $G$ if a subgraph of $G$ is isomorphic to a subdivision of $H$. These two graph containment relations are well-known and well-studied; slightly less ubiquitous is the relation of graph immersion. A pair of adjacent edges $uv,vw$ of a graph $G$ is \textit{lifted} by deleting $uv$ and $vw$ and adding the edge $uw$. If the edge $uw$ already exists in the graph, we add another in parallel to it; if $u=w$ so that the edge $uw$ is a loop, we delete it immediately. We may also lift a path $P$ by deleting all edges of $P$ and adding an edge connecting the ends of $P$. Note that these notions are equivalent -- lifting a pair of edges is a special path-lifting, and each path-lifting is a sequence of pair-liftings. A graph $H$ is called an \textit{immersion} of a graph $G$ if a graph isomorphic to $H$ can be obtained from a subgraph of $G$ by lifting paths. Note that what we have defined here is often called \textit{weak} immersion; all references to the immersion relation in this paper should be understood as referencing weak immersion. From these definitions, we can see that topological minor containment implies both minor and immersion containment, while the minor and immersion containment relations are themselves incomparable. 
	
	While it is true that immersions are relatively unstudied compared to minors and topological minors, it would be inaccurate to state that little has been done. The first mention of immersions was by Nash-Williams \cite{nash2} in 1963, where he conjectured that graphs are well-quasi-ordered by the immersion relation; this was later proven by Robertson and Seymour \cite{RandS2} as part of their Graph Minors project. Several common problems concerning other graph parameters and containment relations have been analogized to immersions, with varying levels of success. For instance, Hadwiger \cite{Hadwiger} famously conjectured that every loopless graph without a $K_t$-minor is $(t-1)$-colorable. The natural analogue for immersions is the conjecture that every graph without a $K_t$-immersion is $(t-1)$-colorable, which was proposed by Lescure and Meyneil \cite{LandM} and later, independently, by Abu-Khzam and Langston \cite{AandL}. This conjecture remains open, with some progress made for specific values of $t$. 
	
	Another problem which has been examined in some detail is that of excluding the containment of certain graphs. In this situation, graph minors are considerably more studied than immersions. In the aforementioned Graph Minors project of Robertson and Seymour \cite{RandS1}, the authors described the general structure of a graph $G$ which does not contain a given graph $H$ as a minor. It is in these papers that the concept of a tree-decomposition of a graph, which we will discuss later, was formalized. In the case of excluding the immersion of a clique, a similar structure theorem was obtained independently by DeVos et al. \cite{DeVos} and Wollan \cite{Wollan}, which utilizes a variation on tree-decompositions based on edge cuts rather than vertex cuts.
	
	Rather than a general graph $H$, one could also try to exclude some set of specific graphs $\{H_1, \ldots, H_n\}$. Perhaps the most famous result of this type is that of Kuratowski \cite{kura}, which states that the class of graphs not containing $K_5$ or $K_{3,3}$ as a topological minor is exactly the class of planar graphs. Numerous results of this nature have been found for minors \cite{ding,oxleyw5} and topological minors \cite{V8}. The number of results on excluding certain immersions, however, is small. DeVos and Malekian provide precise structure theorems for graphs not admitting an immersion of $K_{3,3}$ \cite{nok331} and, in another paper, $W_4$ \cite{noW4}. In her dissertation, Malekian \cite{malekian} also finds results describing graphs not admitting immersions of $K_4$ and the prism graph.
	
	The last problems we would like to discuss, and most relevant to the current paper, concern ``unavoidable substructures". Results of this nature answer questions like ``As a graph $G$ with some connectivity grows large, does some structure necessarily appear in $G$?" These problems are often called Ramsey-type problems, due to their similarity in flavor to the following result of Ramsey \cite{ramsey}.
	
	\begin{thm}
		For every $r\in\mathbb{N}$, there exists an $n\in\mathbb{N}$ such that every simple graph of order at least $n$ contains either $K_r$ or $\overline{K_r}$ as an induced subgraph.
	\end{thm}
	
	In the previous theorem, there is no connectivity requirement on $G$. It is not too difficult to show that if we assume (a simple graph) $G$ to be connected and sufficiently large, then $G$ contains either $K_r, K_{1,r},$ or $P_r$ as an induced subgraph. Such results are very strong in the sense that induced subgraphs are subgraphs, which are themselves topological minors and hence also minors and immersions. Allred, Ding, and Oporowski \cite{sarah} were able to describe the unavoidable induced subgraphs of 2-connected graphs. If we consider $k$-connected graphs for $k\geq 3$, however, the list of unavoidable induced subgraphs quickly grows intractable; for example, any such list for $k=3$ must contain all cubic graphs. It is easier in these settings to relax our containment relation to topological minors, minors, or immersions to keep the list of unavoidable substructures short. Oporowski, Oxley, and Thomas \cite{oot} proved that the unavoidable topological minors of 3-connected graphs are $W_r, L_{r}^+,$ and $K_{3,r}$. In the same paper, they show that the corresponding list for internally 4-connected graphs consists of 5 graphs that they call $A_r, O_r, M_r, K_{4,r}$, and $K_{4,r}'$. The study of unavoidable immersions was essentially nonexistent until Barnes \cite{Matt} showed the following in his dissertation:
	
	\begin{thm}
		For every integer $r\geq 3$, there exists an $n\in\mathbb{N}$ such that every 3-edge-connected graph of order at least $n$ admits an immersion of $L_r^+$ or $P_{2,r}^+$.
	\end{thm}
	
	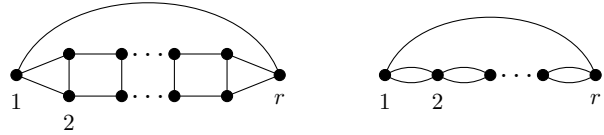
\begin{figure}[ht!]
		\begin{center}
			\begin{tikzpicture}[scale=0.7]
				
				\foreach \x in {0,1,2,3}
				{\draw (\x, 0) -- (\x, .8);
					\foreach \y in {0,.8}
					{\draw[fill=black] (\x,\y) circle (3pt);}}
				
				\foreach \y in {0,.8}
				{\draw (0, \y) -- (1,\y); \draw (2, \y) -- (3,\y);}
				
				\draw (1.5,0) node {$\ldots$};  \draw (1.5,.8) node {$\ldots$};
				
				\draw[fill=black] (-1, .4) circle (3pt);  \draw[fill=black] (4, .4) circle (3pt); 
				\draw (0,.8) -- (-1, .4) -- (0, 0); \draw (3,.8) -- (4, .4) -- (3, 0);
				
				\draw (-1,.4) to[out = 70, in = 110] (4,.4);
				
				\draw (-1,.2) node[below,scale=0.8] {$1$}; \draw (0,-.2) node[below,scale=0.8] {$2$}; \draw (4,.2) node[below,scale=0.8] {$r$};

				\begin{scope}[shift = {(6,.4)}]
					
					\foreach \x in {0,1,2,3,4}{
						\draw[fill=black] (\x, 0) circle (3pt);}
					
					\foreach \x in {0,1,3}{
						\draw (\x,0) to[out = 30, in = 150] (\x+1,0); \draw (\x,0) to[out = -30, in = -150] (\x+1,0); }
					
					\draw (2.5,0) node {$\ldots$};
					\draw (0,0) to[out = 70, in = 110] (4,0);
					
					\draw (0,-.2) node[below,scale=0.8] {$1$}; \draw (1,-.2) node[below,scale=0.8] {$2$}; \draw (4,-.2) node[below,scale=0.8] {$r$};
					
				\end{scope}
			\end{tikzpicture}
			
		\end{center}
		\caption{The graphs $L_r^+$ and $P_{2,r}^+$.}
	\end{figure}
	
	Barnes further conjectured that the only unavoidable immersion of 4-edge-connected graphs was a double cycle of length $r$, denoted $C_{2,r}$. The main result of this paper confirms this conjecture.
	
	\begin{thm}\label{main}
		For every integer $r\geq 4$, there exists a number $f_{\ref{main}}(r)$ such that every $4$-edge-connected graph of order at least $f_{\ref{main}}(r)$ contains $C_{2,r}$ as an immersion.
	\end{thm}
	
	\begin{figure}[ht!]
		\begin{center}
			\begin{tikzpicture}[scale=.6]
				
				\foreach \x in {22.5,67.5, 112.5, 157.5, 202.5, 247.5,337.5}{
					
					\begin{scope}[rotate = \x]
						\coordinate (b) at (67.5:2cm); \coordinate (a) at (112.5:2cm); 
						
						\draw[fill=black] (a) circle (3pt); \draw[fill=black] (b) circle (3pt);
						\draw (a) to[out=30, in =150] (b); \draw (a) to[out=330, in=210] (b);
				\end{scope} }
				
				\dec ([shift=(15:2cm)]0,0) arc (15:30:2cm);
				
			\end{tikzpicture}
		\end{center}
		\caption{The graph $C_{2,r}$.}
	\end{figure}
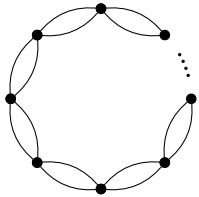
	
	Note that in the previous results, the connectivities of the host graph and the unavoidable graphs are the same; that is, for some small $k$, we consider a $k$-connected host graph and find the unavoidable $k$-connected structures inside of it. However, this quickly becomes difficult for larger values of $k$. This leads to another approach to the study of unavoidable substructures.
	
	Consider the complete bipartite graph $K_{t,r}$ where $t\leq r$. This graph is very well-understood and has many nice properties. It has $t+r$ vertices and $tr$ edges, and it is $t$-connected. What we would like to do is make $K_{t,r}$ our unavoidable substructure; that is, we ask if there exists a function $f(t)$ so that a sufficiently large $f(t)$-connected graph must contain $K_{t,r}$ as a minor. The following result, which follows from the work of B\"{o}hme et al. \cite{bohme}, answers this question in the affirmative.
	
	\begin{thm} \label{bigktk}
		A sufficiently large ($16t+13$)-connected graph contains a $K_{t,r}$ minor.
	\end{thm}
	
	Let $C_{t,r}$ denote the graph obtained from an $r$-vertex cycle by replacing each edge with a parallel class of $t$ edges. This graph is a ``nice" graph for edge-connectivity in the same sense that $K_{t,r}$ is ``nice" for vertex-connectivity. It has $r$ vertices and $tr$ edges, and it is $2t$-edge-connected. Our next result, which we prove in Section 8, is an analogue of Theorem \ref{bigktk}.
	
	\begin{thm}\label{ctr}
		For all positive integers $t,r$, there exists a number $f_{\ref{ctr}}(t,r)$ such that every $(6t-4)$-edge-connected graph of order at least $f_{\ref{ctr}}(t,r)$ contains $C_{t,r}$ as an immersion.
	\end{thm}
	
	Also proven in Section 8 is a rooted version of this theorem, which allows us to choose the terminals of the desired immersion from a specified set $S$ of vertices (provided a slightly higher local edge-connectivity of the set $S$).
	
	In Section 9, we discuss the relationship between immersions, line graphs, and minors; namely, that if $H$ is immersed in $G$, then the line graph of $H$ is a minor of the line graph of $G$. We use this fact, along with Theorem \ref{ctr}, to prove the following.
	
	\begin{thm}\label{linegraph}
		For all positive integers $t,r$, there exists a number $f_{\ref{linegraph}}(t,r)$ such that if $G$ is a $(30t-15)$-connected line graph and $|G| \geq f_{\ref{linegraph}}(t,r)$, then $G$ contains $L(C_{t,r})$ as a minor.
	\end{thm}

	In the forthcoming section, we discuss terminology which we will use throughout this paper. In Section 3, we cite and discuss a result of Mader which is key in our proof. Sections 4 and 5 then discuss tree-decompositions, preparing us for Section 6, where we define a new way to decompose a graph called a ring-decomposition. Section 7 contains the proof of Theorem \ref{main} and Section 8 contains the proof of Theorem \ref{ctr}. Finally, Section 9 contains our discussion on line graphs.

	\section{Terminology}
	
	We begin this section by providing alternate definitions of the topological minor and immersion relations. A graph $H$ is called a \textit{topological minor} of a graph $G$ if some vertices of $G$ can be chosen to represent the vertices of $H$ such that the edges of $H$ can then be represented by internally vertex-disjoint paths of $G$. A graph $H$ is called an \textit{immersion} of a graph $G$ if some vertices of $G$ can be chosen to represent the vertices of $H$ such that the edges of $H$ can then be represented by edge-disjoint paths of $G$. Note the similarity between these definitions; the changing of internally vertex-disjoint paths to edge-disjoint paths makes the immersion relation a natural place to study edge-analogues of common graph theory problems. We will call the vertices of $G$ chosen to represent the vertices of $H$ the \textit{terminals} of the immersion. 
	
	Next we will define the notion of \textit{bridges} of a subgraph. Note that by the term bridges we do not mean, as is sometimes the case, cut-edges. Instead, we use the notion set out by Tutte \cite{tuttebridges}; for this reason, the subgraphs we will define are sometimes known as \textit{Tutte bridges}. Suppose that $H$ is a subgraph of $G$. An \textit{$H$-bridge} in $G$ is a connected subgraph $B$ of $G$ such that either $B$ is a single edge of $E(G)\bs E(H)$ with both endpoints contained in $V(H)$, or $B-V(H)$ is a connected component of $G-V(H)$ and $B$ contains every edge of $G$ with one end in $V(B)\bs V(H)$ and the other in $V(H)$.
	In the case that $B$ is a single edge with both ends in $H$, we call $B$ a \textit{trivial bridge}. In either case, the vertices of $V(B)\cap V(H)$ are called the \textit{attachments} of $B$ and denoted $Att(B)$. The graph $B-Att(B)$ is called the \textit{nucleus} of $B$. We will use the notion of bridges in the proof of Lemma \ref{claim}.
	
	We now define several graphs we will reference throughout the paper. Let $K_{t,r}$ denote the complete bipartite graph with partite sets of size $t$ and $r$. Let $C_r$ denote the cycle on $r$ vertices, and let $C_{t,r}$ denote the graph obtained by replacing each edge of $C_r$ by a parallel class of edges of size $t$. Similarly, let $P_r$ denote the path on $r$ vertices and $P_{t,r}$ the graph obtained by replacing each edge of $P_r$ by a parallel class of edges of size $t$.

	Now suppose that $A$ and $B$ are sets. Then by $A\bs B$ we mean the set of elements which are in $A$ and are not in $B$. If instead we suppose that $A$ is a graph, then \begin{itemize}
		\item if $B\subseteq V(A)$, then $A - B$ denotes the graph obtained from $A$ by deleting the vertices in $B$;
		\item if $B\subseteq E(A)$, then $A\bs B$ denotes the graph obtained from $A$ by deleting the edges in $B$;
		\item if $B$ is a subgraph of $A$, then $A-B$ denotes $A - V(B)$ and $A\bs B$ denotes $A \bs E(B)$.
	\end{itemize}
	
	Next, for a set of vertices $X$ in a graph $G$, we denote by $N(X)$ the set of neighbors of $X$, that is, the set of vertices not in $X$ that are adjacent to some member of $X$. Similarly, let $\delta(X)$ denote the set of edges between a member of $X$ and a non-member of $X$. In a slight abuse of notation, we often write refer to set $\{x\}$ as simply $x$.
	
	Lastly, for a positive integer $n$, we let $[n]$ denote the set $\{1,\ldots, n\}$. Any notation not explicitly detailed here will follow the conventions of Diestel's \cite{diestel} text.

	\section{Maximum degree}
	
	In researching unavoidable substructures, it is common practice to simplify the given graph as much as possible, allowing underlying structures to be seen more easily. In this section, we will show how an important result of Mader lets us bound the maximum degree of graphs.
	
	Let $G = (V, E)$ be a graph with $u,v\in V$ and let $\lambda(u,v; G)$ denote the local edge-connectivity between $u$ and $v$, that is, the maximum number of edge-disjoint paths connecting $u$ and $v$ in $G$. Let $e = su$ and $f=sv$ be two distinct edges of $G$. Recall that the operation of deleting $e$ and $f$ and adding the edge $uv$ is known as lifting the pair of edges $\{e,f\}$. Let $G^{ef}$ denote the graph obtained by lifting the pair $\{e,f\}$. If $\lambda(x,y; G^{ef}) = \lambda(x,y; G)$ for all $x,y\in V\bs \{s\}$, then the pair of edges $\{e,f\}$ is called \textit{liftable}.
	
	With these definitions in hand, we are prepared to state the following result of Mader \cite{mader}.
	
	\begin{thm}\label{mader}
		Let $G= (V, E)$ be a connected graph such that $s\in V$ is not incident to a cut-edge of $G$ and $d(s)\neq 3$. Then there is a liftable pair $\{e,f\}$ of edges incident to $s$.
	\end{thm}
	
	The following observation will be key to simplifying the proof of our main result as well as many other arguments concerning immersions: 
	
	\begin{obs} \label{obs1}
		Suppose $G$ is a $k$-edge-connected graph where $k>1$. Then $G$ contains an immersion $H$ such that $H$ is $k$-edge-connected, every vertex of $H$ has degree $k$ or $k+1$, and $V(H)=V(G)$.
	\end{obs}
	
	\noindent\textit{Proof:} Note that $G$ has minimum degree at least $k$. By fixing a vertex $u$ of degree at least $k+2$, Mader's theorem implies that we can lift a pair of edges adjacent to $u$ to produce a graph $G'$ where $d_{G'}(u) = d_G(u) - 2 \geq k$, $V(G') = V(G)$, and $\lambda(x,y; G') = \lambda(x,y; G) \geq k$ for all vertices $x,y$ not equal to $u$. However, we must ensure that $\lambda(u,v;G') \geq k$ for all $v\neq u$.
	
	Suppose that in $G'$ this is not the case. Then $G'$ contains an edge-cut $X$ of size less than $k$ between $u$ and $v$ for some $v\in V(G')\bs u$. But since Mader's theorem guarantees $\lambda(w,v;G') = \lambda(w,v;G) \geq k$ for $w,v\in V(G')\bs u$, we must have that every such pair are in the same component of $G'\bs X$. But this implies that $X = \delta(u)$, contradicting the fact that the degree of $u$ in $G'$ is at least $k$. Therefore we have that $V(G') = V(G)$ and $G'$ is $k$-edge-connected, and hence we may repeatedly apply Mader's theorem to each vertex of $G$ in turn until the desired immersion $H$ is reached.  \hfill $\Box$ \\
	
	Consider a large 4-edge-connected graph $G$. Because of the previous observation, we know that $G$ has an immersion $H$ on the same vertex set which is 4-edge-connected and has only vertices of degree 4 or 5. If we can prove that $H$ contains the desired $C_{2,r}$ immersion, then the transitivity of the immersion relation implies that our original graph $G$ does as well. In our application, therefore, we may assume that every vertex of the graph $G$ has degree 4 or 5.

	\section{tree-decompositions}
	
	Informally, a tree-decomposition of a graph $G$ is a tool which yields some measure (called tree-width) of how similar $G$ is to a tree. To be precise, let $G$ be a graph, $T$ be a tree, and $\mathcal{Y} = \{Y_t\}_{t\in V(T)}$ be a family of subsets of $V(G)$ indexed by the vertices of $T$. Then the pair $(T,\mathcal{Y})$ is called a \textit{tree-decomposition} of $G$ if: 
	
	\begin{enumerate}[label=(W{\arabic*})]
		
		\item $\bigcup\limits_{t\in V(T)} Y_t = V(G)$, and every edge of $G$ has both ends in some $Y_t$; and
		
		\item if $t,t',t''\in V(T)$ and $t'$ lies on the path between $t$ and $t''$, then $Y_t\cap Y_{t''} \subseteq Y_{t'}$.
		
	\end{enumerate}
	
	The \textit{width} of a tree-decomposition $(T, \mathcal{Y})$ is $\text{max}_{t\in V(T)}(|Y_t| - 1)$, and the \textit{tree-width} of $G$ is the minimum width over all tree-decompositions of $G$. The members of $\mathcal{Y}$ are called the \textit{bags} of the tree-decomposition.
	
	The following theorem of Chudnovsky et. al \cite{chudnov} provides the first instance in which we can guarantee the existence of the immersion we desire.
	
	\begin{thm}\label{chudnovsky}
		For all $g > 1$ there exists a number $f_{\ref{chudnovsky}}(g)$ such that every $4$-edge-connected graph with tree-width at least $f_{\ref{chudnovsky}}(g)$ contains $J_g$ as an immersion.
	\end{thm}
	
	Here, $J_g$ denotes the $g\times g$ grid graph. We state the following corollary for use in our specific application.
	
	\begin{coro}\label{chud}
		For any integer $r\geq 2$, there exists a number $f_{\ref{chud}}(r)$ such that every $4$-edge-connected graph with tree-width at least $f_{\ref{chud}}(r)$ contains $C_{2,r}$ as an immersion.
	\end{coro}
	
\begin{figure}[ht!]
		\begin{center}
			\begin{tikzpicture}[scale=0.5]
				
				\foreach \x in {0,1,2,3,4,5,6}{
					\draw[opacity=0.5] (\x, 0) -- (\x, 6); \draw[opacity=0.5] (0, \x) -- (6,\x);}

				\draw[very thick] (1,5) -- (1,6) -- (6,6) -- (6,1) -- (5,1);
				\draw[very thick] (1,5) -- (0,5) -- (0,0) -- (5,0) -- (5,1);
				
				\draw[very thick] (1,5) -- (2,5) -- (2,4) -- (1,4) -- (1,5);
				\draw[very thick] (2,4) -- (3,4) -- (3,3) -- (2,3) -- (2,4);
				\draw[very thick] (3,3) -- (4,3) -- (4,2) -- (3,2) -- (3,3);
				\draw[very thick] (4,2) -- (5,2) -- (5,1) -- (4,1) -- (4,2);
				
				\draw (1,5) node[above right,scale=0.75] {$v_1$}; 
				\draw (2,4) node[above right,scale=0.75] {$v_2$};
				\draw (3,3) node[above right,scale=0.75] {$v_3$};
				\draw (4,2) node[above right,scale=0.75] {$v_4$};
				\draw (5,1) node[above right,scale=0.75] {$v_5$};
				
				\foreach \x in {0,1,2,3,4,5,6}
				{\foreach \y in {0,1,2,3,4,5,6}
					\draw[fill=black,opacity=0.5] (\x,\y) circle (1.5pt);}
				
				\draw[fill=black] (1,5) circle (3pt); \draw[fill=black] (2,4) circle (3pt);
				\draw[fill=black] (3,3) circle (3pt); \draw[fill=black] (4,2) circle (3pt);
				\draw[fill=black] (5,1) circle (3pt);
				
			\end{tikzpicture}
		\end{center}
		\caption{A $C_{2,5}$ immersion in $J_{7}$.}
		\label{grid}
	\end{figure}
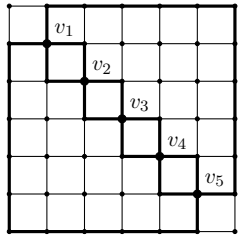

	\noindent\textit{Proof:} Let $f_{\ref{chud}}(r) := f_{\ref{chudnovsky}}(r+2)$. The result follows because $J_{r+2}$ contains a $C_{2,r}$ immersion (see Figure \ref{grid}).\\
	
	\section{More conditions on a tree-decomposition}
	
	Additional conditions can be imposed on a tree-decomposition $(T,\mathcal{Y})$ to reduce redundancy and to reveal more structures. We will consider the following properties, which were introduced in \cite{oot}.
	
	\begin{enumerate}[label=(W{\arabic*})]
		\setcounter{enumi}{2}
		
		\item for every two vertices $t, t'$ of $T$ and every positive integer $k$, either there are $k$ disjoint paths in $G$ between $Y_t$ and $Y_{t'}$, or there is a vertex $t''$ of $T$ on the path between $t$ and $t'$ such that $|Y_{t''}| < k$.
		
		\item if $t, t'$ are distinct vertices of $T$, then $Y_t \neq Y_{t'}$.
		
		\item if $t_0\in V(T)$ and $B$ is a component of $T-t_0$, then $\bigcup\limits_{t\in V(B)} Y_t \bs Y_{t_0} \neq \emptyset.$
		
	\end{enumerate} 
	
	Let $G$ be a connected graph with $\Delta(G) \leq d$ and let $(T,\mathcal{Y})$ be a tree-decomposition of $G$ of width at most $w$ which satisfies (W5). Then:
	
	\begin{remark}\label{w2w5}
		$\Delta(T) \leq (w+1)d$.
	\end{remark}
	
	\noindent \textit{Proof:} Consider an arbitrary vertex $t_0$ of $T$, and suppose that the degree of $t_0$ in $T$ is $s$. Then the removal of $t_0$ from $T$ yields $s$ components $T_1, T_2, \ldots, T_s$ of $T- t_0$. For each component $T_i$, define $B_i = \bigcup\limits_{t\in T_i} Y_t \bs Y_{t_0}.$ Then (W5) implies that $B_i$ is not empty for each $i\in[s]$. Furthermore, (W2) implies that the $B_i$ are pairwise-disjoint, since any vertex in the intersection of a $B_i$ and a $B_j$ where $i\neq j$ would necessarily belong to $Y_{t_0}$. Note that this implies that $Y_{t_0}, B_1, \ldots, B_s$ form a partition of $V(G)$.
	
	Now we argue that no edge of $G$ exists between the subsets $B_i$ for $i\in[s]$. Suppose there exists an edge $e=xy$ where $x\in B_m$ and $y\in B_n$ for $m\neq n$. Then both ends of $e$ must be contained in some bag $Y_t$ for some $t$ (by (W1)), contradicting the fact that $Y_{t_0}, B_1, \ldots, B_s$ partitions $V(G)$. Hence no such edge $e$ exists. 
	
	Lastly, we note that since $G$ is a connected graph, it must contain an edge leaving each $B_i$ (that is, $\delta(B_i)$ is nonempty for each $i\in[s]$). But the argument above shows that any such edge must be between $B_i$ and $Y_{t_0}$. This implies that each $B_i$ may be uniquely associated with a member of $\delta(Y_{t_0})$. However, since there are at most $(w+1)$ vertices in $Y_{t_0}$ and each has degree at most $d$, we find that $|\delta(Y_{t_0})| \leq (w+1)d$. This implies that there are at most $(w+1)d$ of the $B_i$, and hence that $s\leq (w+1)d$, completing the proof. \hfill $\Box$ \\
	
	The arguments of the previous proof also imply the following, which we note here and make use of later.
	
	\begin{remark}\label{components}
		For any $t\in T$, $G - Y_t$ has at least $d_T(t)$ components.
	\end{remark}
	
	Tree-decompositions are often called \textit{linked} \cite{linked} when they satisfy (W3). In \cite{oot}, the authors call a tree-decomposition \textit{lean} when it satisfies (W1)-(W5) as well as an additional property (W6); since we have no need for this final property, we shall refer to a tree-decomposition as \textit{lean} when it satisfies (W1)-(W5). 
	
	The following result of Oporowski, Oxley, and Thomas \cite{oot} allows us to consider lean tree-decompositions in our proof:
	
	\begin{thm} \label{lean}
		If a graph has treewidth $w$, then that graph has a lean tree-decomposition of width at most $w$.
	\end{thm}
	
	Now suppose that the tree-decomposition $(T,\mathcal{Y})$ is lean and that $T$ contains a long path $P$; we can show that $P$ contains a long sequence of bags with special properties. Note that the first part of the following assertion and its proof are simply Claim (4) from Theorem 3.1 of \cite{oot}.
	
	\begin{lemma}\label{railsgates}
		There exists a function $f_{\ref{railsgates}}(n,w)$ such that, for all positive integers $n$ and $w$, the following holds: If $G$ is a connected graph and $(T, \mathcal{Y})$ is a lean tree-decomposition of $G$ having width at most $w$ and containing a path $P$ of length at least $f_{\ref{railsgates}}(n,w)$, then $V(P)$ contains a subset $\{t_1, t_2, \ldots, t_{n}\}$ of interior vertices occurring along $P$ in the order listed such that:
		
		\begin{enumerate}[i.]
			\item for some positive integer $s \leq w+1$, $|Y_{t_i}| = s$ for all $i\in[n]$ and $|Y_t| \geq s$ for every $t$ in $V(P)$ between $t_1$ and $t_n$; and 
			\item $Y_{t_i}\cap Y_{t_j} = U$ for some $U\subseteq V(G)$ and for all distinct $i,j\in [n]$.
		\end{enumerate}

	\end{lemma}
	
	\noindent\textit{Proof:} Let $f_{\ref{railsgates}}(n,w) := 2+n_1^{w+1}$ where $n_1 = n^{w+1}$. Suppose that $G$ is a connected graph and $(T, \mathcal{Y})$ is a lean tree-decomposition of $G$ having width at most $w$ and containing a path $P$ of length at least $f_{\ref{railsgates}}(n,w)$. We first claim that $P$ contains a sequence of interior vertices of length at least $n_1$ which satisfies (i). Because $G$ is connected and each interior vertex of $P$ represents a vertex cut of $G$ (as noted in Remark \ref{components}), we must have that $|Y_t|\geq 1$ for all interior vertices $t$ of $P$. Now we find $n_1$ disjoint subpaths of $P$ each containing exactly $(n_1)^w$ interior vertices of $P$. In each such subpath, we look for vertices $t$ such that $|Y_t| = 1$; if we find one such bag in each subpath, we have satisfied the claim. In the case that we do not find such a bag in some subpath $P'$, we must have that $|Y_t| \geq 2$ for all $t\in V(P')$. We then find $n_1$ disjoint subpaths of $P'$ each containing $(n_1)^{w-1}$ vertices, and we search these for vertices $t$ such that $|Y_t| = 2$. Continuing to argue in this way, we conclude that, since $|Y_t|\leq w+1$ for all $t\in V(P)$, there must be some integer $s$ for which our claim holds. Moreover, we see that $1\leq s \leq w+1$.
	
	We now claim that any sequence of bags of length at least $n^s$ which satisfies (i) contains a subsequence of bags of length $n$ which satisfies (ii). We prove this by induction on $s$. If $s=1$, then property (W4) of $(T, \mathcal{Y})$ implies that (ii) holds (with $U = \emptyset$) for any sequence of $n$ bags. Now suppose that $s>1$. If some vertex $v$ of $G$ belongs to at least $n^{s-1}$ of the $Y_{t_i}$, then the result follows upon removing $v$ from each of the $Y_{t_i}$, applying the inductive hypothesis to find the correct subsequence, and then replacing $v$. If each vertex of $G$ belongs to fewer than $n^{s-1}$ of the $Y_{t_i}$, then $t_1, t_{n^{s-1}}, t_{2n^{s-1}},\ldots, t_{(n-1)n^{s-1}}$ yields the desired subset of vertices of $P$ (in particular, they are disjoint by (W2)). Lastly, we note that $n^s \leq n^{w+1} = n_1$, and hence the argument in the previous paragraph implies that $P$ contains a sequence of interior vertices $t_1, \ldots, t_{n^s}$ which satisfies (i). We may then apply our claim to finish the proof. \hfill $\Box$ \\
	
	Note that we may actually take $s< w+1$ in Lemma \ref{railsgates}, a fact that follows from the leanness of the tree-decomposition. However, to prove this requires a slightly more complicated argument, and we opt for the present version since it suffices for our applications.
	
	Moving forward, we will refer to the specified bags of size $s$ in Lemma \ref{railsgates} as \textit{gates} of $(G, (T, \mathcal{Y}))$ or simply as gates of $G$ when $(T,\mathcal{Y})$ is fixed. Let $P'$ be the subpath of $P$ with ends $t_2$ and $t_{n-1}$. Suppose that $(T,\mathcal{Y})$ is lean (Theorem \ref{lean} allows us to assume that this is the case). Then $G$ has $s$ vertex-disjoint paths $R_1,R_2, \ldots, R_s$ between $Y_{t_1}$ and $Y_{t_n}$. Fix $s$ and these disjoint paths. If $t,t'\in V(P')$ are gates, then for $j\in[s]$ there is a unique subpath of $R_j$ which has one end in $Y_t$ and the other in $Y_t'$. We denote this subpath by $R_j(t,t')$. Observe that if $v_1, v_2, \ldots, v_p$ lie on $P'$ in this order and $|Y_{v_i}| = s$ for $i \in[p]$, then $R_j(v_1,v_p)$ is obtained by pasting together $R_j(v_1,v_2), R_j(v_2, v_3), \ldots, R_j(v_{p-1}, v_p)$ in this order. We will refer to the paths $R_j$ as \textit{rails} of $(G, (T, \mathcal{Y}))$ or simply as rails of $G$ when $(T,\mathcal{Y})$ is fixed. \\
	
	\section{Ring-decompositions}
	
	Let $G$ be a graph and $w$ be a positive integer. A \textit{ring-decomposition} of $G$ of \textit{width} $w$ and \textit{length} $n$ is a cyclic sequence $(G_0, G_1, \ldots, G_n)$ ($n\geq 3$) of subgraphs of $G$ such that:
	\begin{enumerate}[label=(R{\arabic*})]
		\item $G_0 \cup G_1 \cup \ldots \cup G_n = G$
		\item consecutive $G_i$ are edge-disjoint and non-consecutive $G_i$ are vertex-disjoint
		\item for each $i\in[n+1],$ let $W_i := V(G_{i-1}\cap G_i)$, where $G_{n+1} = G_0$; then $|W_i| = w$
		\item for each $i\in[n],$ $G_i$ contains $w$ disjoint paths between $W_i$ and $W_{i+1}$.
	\end{enumerate}
	
	Note that (R4) need not apply to $G_0$, which differentiates it from $G_1,\ldots, G_n$. When working with ring-decompositions, we often would like to impose some further properties on the subgraphs $G_1, \ldots, G_n$. To this end, we say a ring-decomposition $(G_0, G_1,\ldots, G_n)$ is \textit{connected} if the subgraph $G_i$ is connected for every $i\in[n]$ (note that $G_0$ may not be connected). The following structural theorem is the main result of this section.
	
	\begin{thm}\label{connected}
		There exists a function $f_{\ref{connected}}(n,w,d)$ satisfying the following. For any integer $n\geq 3$, if $G$ is a connected graph with $tw(G)\leq w$, $\Delta(G) \leq d$, and $|G| \geq f_{\ref{connected}}(n,w,d)$, then $G$ has a connected ring-decomposition of width at most $w$ and length at least $n$.
	\end{thm}
	
	Before proving this result, we note one way to simplify our arguments. By taking the union of some subset of consecutive $G_i$, we can create a ring-decomposition of the same width having shorter length. We omit the straightforward proof of the following remark.
	
	\begin{remark}\label{absorb}
		Let $\mathcal{G} = (G_0, G_1,\ldots, G_n)$ be a (connected) ring-decomposition of a graph $G$ having width $w$. Then by letting $G_{n+1} = G_0$, we have that for any $i\in[n]$, the sequence of subgraphs formed by replacing the two terms $G_i, G_{i+1}$ of $\mathcal{G}$ with $G_i\cup G_{i+1}$ forms a (connected) ring-decomposition of width $w$ and length $n-1$, provided $n-1\geq 3$.
	\end{remark}
	
	In practice, we can apply this operation many times to find that any consecutive subsequences of our $G_i$ can be combined to form a (possibly much) shorter ring-decomposition with some useful property. We will use this fact several times to simplify the proof of Theorem \ref{connected}, which we reach through a series of structural lemmas. The first such lemma is the following, which is implied by Theorem 5.3 of \cite{binary}.
	
	\begin{lemma} \label{paththm}
		For every positive integer $r$, there is a number $f_{\ref{paththm}}(r)$ such that every connected graph of order at least $f_{\ref{paththm}}(r)$ contains $K_{1,r}$ or $P_r$ as a subgraph.
	\end{lemma}
	
	In the following three lemmas, let $G$ be a connected graph with $\Delta(G) \leq d$ and let $(T,\mathcal{Y})$ be a lean tree-decomposition of $G$ of width at most $w$.
	
	\begin{lemma} \label{longpath}
		There exists a function $f_{\ref{longpath}}(n,w,d)$ such that, for all positive integers $n,w,d$, the following holds: If $|G| \geq f_{\ref{longpath}}(n,w,d)$, then $T$ contains a path $P$ and a sequence $\{t_1,\ldots, t_n\}$ of interior vertices of $P$ satisfying (i)--(ii) of Lemma \ref{railsgates}.
	\end{lemma}
	
	\noindent\textit{Proof:} Remark \ref{w2w5} shows that $\Delta(T) \leq (w+1)d$. Let $n' = \max\{f_{\ref{railsgates}}(n,w),(w+1)d+1\}$. Let $f_{\ref{longpath}}(n,w,d) := (f_{\ref{paththm}}(n'))(w+1)$. Then 
	\begin{align*}
		|V(T)|\cdot(w+1) &\geq |V(G)|, \text{ and hence} \\
		|V(T)|\cdot (w+1) &\geq (f_{\ref{paththm}}(n'))(w+1), \text{ so that} \\
		|V(T)| &\geq f_{\ref{paththm}}(n').
	\end{align*}
	
	Then by applying Lemma \ref{paththm} to $T$, we can guarantee that $T$ contains either $K_{1,n'}$ or $P_{n'}$ as a subgraph. However, since $\Delta(T) < n'$, we cannot have that $T$ contains $K_{1,n'}$. Therefore $T$ must contain $P_{n'}$ and hence $P_{f_{\ref{railsgates}}(n,w)}$ as a subgraph. Lemma \ref{railsgates} then implies the result. \hfill $\Box$ \\
	
\begin{figure}[ht!] 
		\begin{center}
			\begin{tikzpicture}[scale=0.5]
				
				\draw[thick] (-2,0) -- (10,0);
				
				\draw[dashed] plot [smooth cycle] coordinates {(-.5,-1) (8.25,-1) (8.25,.5) (7,1) (6.75,2.5) (-1,2.5) (-.75,.5)};
				
				\draw (0,-.5) node {$t'$}; \draw (8,-.5) node {$t''$};
				\draw (-1.75, .5) node {$P$};

				\draw (-.75, 2) -- (-.5,1)--(0,0)--(.5,1)--(.75,2); \draw (-.5,1)--(-.25,2); \draw (.25,1.5)--(.5,1);
				
				\draw (1.75,1.5)--(2,1)--(2,0);
				
				\draw (5.5,2)--(6,1)--(6,0); \draw (5.5,1)--(6,0); \draw (6.5,2)--(6,1);
				
				\draw (8,2)--(8.5,1)--(8,0); \draw (9,2)--(8.5,1);
				
				\foreach \x in {-.75,-.25,.75,5.5,6.5,8,9}{
					\draw[fill=black] (\x,2) circle (3pt);}
				
				\foreach \x in {.25,1.75}{
					\draw[fill=black] (\x,1.5) circle (3pt);}
				
				\foreach \x in {-.5,.5,2,5.5,6,8.5}{
					\draw[fill=black] (\x,1) circle (3pt);}
				
				\foreach \x in {0,2,4,6,8}{
					\draw[fill=black] (\x,0) circle (3pt);}

			\end{tikzpicture}
		\end{center}
		\caption{The set $Y[t',t'')$ is the union of the bags within the dotted outline.}
		\label{bags}
	\end{figure}
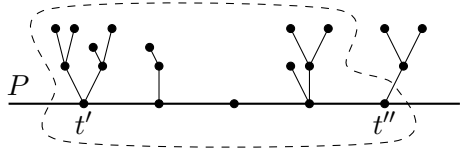

	In the next two lemmas, let $P$ and $\{t_i\}$ be as described in Lemma \ref{longpath}. Let $P[t',t'']$ denote the subpath of $P$ with ends $t'$ and $t''$. Then we define $Y_P[t', t'')$ (see Figure \ref{bags}) as the union of all bags $Y_t$ such that:
	\begin{itemize}
		\item $t\in P[t',t'']$, or
		\item $t\in V(T-P)$ and the unique path in $T$ from $t$ to $P$ ends in a vertex of $P[t',t'']-t''$.
	\end{itemize}
	
	When $P$ is fixed or the omission will otherwise not cause confusion, we refer to $Y_P[t',t'')$ simply as $Y[t',t'')$.
	
	\begin{lemma}\label{disjointY}
		For any $1\leq i < i' < j < j' \leq n$, the sets $Y[t_i, t_{i'}) \bs U$ and $Y[t_j, t_{j'})\bs U$ are disjoint.
	\end{lemma}
	
	\noindent\textit{Proof:} Suppose that $1\leq i < i' < j < j' \leq n$ and there exists some vertex $v\in (Y[t_i, t_{i'}) \bs U)\cap (Y[t_j, t_{j'})\bs U)$. Then, since $v$ induces a connected subgraph of $T$ (property (W2) of tree-decompositions), we have that $v\in Y_{t_{i'}}$ and that $v\in Y_{t_j}$. But this means $v\in Y_{t_{i'}}\cap Y_{t_j}$; since $i'\neq j$, we have that $Y_{t_{i'}}\cap Y_{t_j} = U$, a contradiction. \hfill $\Box$ \\

	\begin{lemma}\label{noU}
		
		For every integer $m>1$, there is a number $\gamma = f_{\ref{noU}}(m,w,d)$ such that any sequence of gates of length at least $\gamma$ which satisfies (i)--(ii) of Lemma \ref{railsgates} contains, for some $k$, a subsequence of gates $\{Y_{t_{k+1}}, \ldots, Y_{t_{k+m}}\}$ such that $N(U) \cap (Y[t_{k+1},t_{k+m})\bs U) = \emptyset$.
	\end{lemma}
	
	\noindent\textit{Proof:} Let $\gamma = f_{\ref{noU}}(m,w,d) = (dw+1)m$. Note that because the width of $(T,\mathcal{Y})$ is at most $w$, we have that $|U| \leq w+1$. Note that $d,w$ and $m$ are positive integers, and hence $\gamma > 1$; because $(T,\mathcal{Y})$ satisfies (W4), we see that $|U| \leq w$. This fact, along with the assumption that $\Delta(G)\leq d,$ implies that $|N(U)| \leq dw$. Now we may divide the first $(dw+1)m$ gates into $dw+1$ consecutive subsequences of length $m$. Consider any two such subsequences $\{Y_{t_{im+1}}, \ldots, Y_{t_{(i+1)m}}\}$ and $\{Y_{t_{jm+1}},\ldots, Y_{t_{(j+1)m}}\}$. Lemma \ref{disjointY} implies that $Y[t_{im+1}, t_{(i+1)m})\bs U$ and $Y[t_{jm+1},t_{(j+1)m})\bs U$ are disjoint. Finally, since there are at most $dw$ neighbors of the set $U$ and we have $dw+1$ disjoint consecutive subsequences of length $m$, at least one must contain no neighbors of $U$, as desired. \hfill $\Box$ \\
	
	We are now prepared to show that a sufficiently large connected graph with bounded degree and bounded tree-width admits a ring-decomposition. 
	
	\begin{lemma}\label{ringdecomp}
		
		There exists a function $f_{\ref{ringdecomp}}(n,w,d)$ satisfying the following. For any integer $n\geq 3$, if $G$ is a connected graph with $tw(G)\leq w$, $\Delta(G) \leq d$, and $|G| \geq f_{\ref{ringdecomp}}(n,w,d)$, then $G$ has a ring-decomposition of width at most $w$ and length at least $n$.
		
	\end{lemma}
	
	\noindent\textit{Proof:} Suppose that $G$ is a connected graph with $tw(G) \leq w$ and $\Delta(G) \leq d$. Let $n_1 = f_{\ref{noU}}(n+1,w,d)$. We claim that if $|G| \geq f_{\ref{ringdecomp}}(n,w,d) : = f_{\ref{longpath}}(n_1, w, d)$, then $G$ has the desired ring-decomposition. 
	
	Firstly, Lemmas \ref{lean} and \ref{longpath} tell us that $G$ has a lean tree-decomposition $(T,\mathcal{Y})$ of width at most $w$ which contains a sequence of $n_1$ gates which satisfies (i)--(ii) of Lemma \ref{railsgates}; these gates have common size $s$, where $s\leq w+1$. Lemma \ref{noU} then implies that this sequence of gates has a subsequence $\{Y_{t_1}, \ldots, Y_{t_{n+1}}\}$ with intersection $U$ such that $N(U) \cap (Y[t_1,t_{n+1})\bs U) = \emptyset$.
	
	For $t,t'\in\{t_1,\ldots,t_{n+1}\}$, define $G[t,t'] := G[Y[t,t')\bs U]\bs E(G[Y_{t'}])$. Let $G_0'$ be the subgraph of $G$ induced by the set of edges $E(G)\bs E(G[t_1, t_{n+1}])$, and define $G_0$ as the subgraph of $G$ having vertex set $V(G_0')\cup Y_1 \cup Y_{t_{n+1}}$ and edge set $E(G_0')$. We claim that the graphs $G_i := G[t_i, t_{i+1}]$ for $i\in[n]$ along with $G_0$ form a ring-decomposition of $G$ of width $s-|U|$ and length $n$.
	
	To show this, we need only show that $(G_0, G_1, \ldots, G_n)$ satisfy properties (R1) - (R4) stated in the definition of ring-decompositions. It is clear from the definition of $G_0$ that every edge of $G$, and hence every vertex of $G$, is in some $G_i$ or $G_0$. Thus $G_0 \cup G_1 \cup \ldots \cup G_n = G$ and (R1) is satisfied. 
	
	It is clear from the definition of $G_0$ that $G_0$ shares vertices only with $G_1$ and $G_n$. Recall that for $i\in[n]$, we have $G_i = G[t_i, t_{i+1}] = G[Y[t_i, t_{i+1}) \bs U]\bs E(G[Y_{t_{i+1}}])$. This means that $V(G_i) = Y[t_i, t_{i+1}) \bs U$, and hence Lemma \ref{disjointY} implies that non-consecutive $G_i$ are vertex-disjoint.
	
	Now consider consecutive subgraphs $G_i,G_{i+1}$. We see $V(G_i) = Y[t_i, t_{i+1})\bs U$ and $V(G_{i+1}) = Y[t_{i+1}, t_{i+2})\bs U$, which indicates that $Y_{t_{i+1}} \bs U$ is contained in the vertex set of both subgraphs. Suppose that for some vertex $v\notin Y_{t_{i+1}} \bs U$ we have $v\in G_i\cap G_{i+1}$. Then $v\in Y^a$ for some bag $Y^a$ contained in $T_1$, the subtree of $T$ corresponding to $G_i$; that is, the subgraph of $T$ induced by the set of bags $B$ forming the union $Y[t_i, t_{i+1})$. Similarly, we have that $v\in Y^b$ for some $Y^b$ contained in $T_2$, the subtree of $T$ corresponding to $G_{i+1}$. However, property (W2) of tree-decompositions then implies that $Y^a$ and $Y^b$ are both contained in $T_v$, the subtree of $T$ induced by all bags containing $v$. But this further implies that $v\in Y_{t_{i+1}}$, a contradiction  (note that every vertex of $U$ is in $G_0$). Therefore we find that $V(G_i)\cap V(G_{i+1}) = Y_{t_{i+1}} \bs U$ for all $i\in[n-1]$. Note that the definition of $G_0$ immediately yields the facts $V(G_0) \cap V(G_1) = Y_{t_1}$ and $V(G_0)\cap V(G_n) = Y_{t_{n+1}}$, confirming that the decomposition $(G_0, G_1, \ldots, G_n)$ satisfies (R3) with width $s - |U|$. To finish establishing that (R2) holds, we need only confirm that $G_i$ and $G_{i+1}$ are edge-disjoint. This follows easily from the fact that they intersect only in the vertices $Y_{t_{i+1}} \bs U$ and the edges among this set belong (by definition) to $G_{i+1}$ and not $G_i$. 
	
	Lastly, we note that the existence of the rails purported after Lemma \ref{railsgates} implies property (R4). \hfill $\Box$ \\
	
	It is now relatively simple to prove Theorem \ref{connected}: \\
	
	\noindent\textit{Proof of Theorem \ref{connected}}: Suppose that $G$ is a connected graph with $tw(G) \leq w, \Delta(G) \leq d,$ and $|G| \geq f_{\ref{connected}}(n,w,d) := f_{\ref{ringdecomp}}(n^{2^p}, w, d)$, where $p = \lceil \log_2 w\rceil$.
	
	We prove this theorem by way of a simpler statement. Suppose that a connected graph $G$ has a ring-decomposition $G_0, G_1, \ldots, G_{m^2}$ of width $\omega$. Then either $G$ has a connected ring-decomposition of width $\omega$ and length $m$, or $G$ admits a ring-decomposition having width at most $\omega/2 $ and length $m$. 
	
	To see this, simply cut the sequence $G_1, \ldots, G_{m^2}$ into $m$ disjoint subsequences of length $m$. Next, consider the union of each subsequence, say $H_1 := G_1 \cup \ldots \cup G_m$, $H_2 := G_{m+1} \cup \ldots \cup G_{2m}$, etc. If every graph $H_i$ is connected, then $(G_0, H_1, \ldots, H_m)$ is a connected ring-decomposition of $G$ having width $\omega$ (as noted in Remark \ref{absorb}). If this is not the case, then some $H_i$ is disconnected; suppose for simplicity that this is $H_1$ (the argument for other $H_i$ is analogous). Let $R'$ denote the set of the restrictions of the rails of $G$ to $H_1$. Let $K$ denote the component of $H_1$ containing the fewest members of $R'$, and let the set of these members of $R'$ contained in $K$ be denoted $R$. Since the members of $R'$ are partitioned among the components of $H_1$, we must have that $|R| \leq \omega/2$. For $i\in[m]$, let $G_i' = G_i \cap K$. Then $K = G_1'\cup \ldots \cup G_m'$ and we define $G_0'$ to be the subgraph of $G$ induced by all edges of $G$ which are not in $K$. It is then straightforward to verify that $(G_0', G_1', \ldots, G_m')$ forms a ring-decomposition of width $|R|$, which we have noted is at most $\omega/2$.
	
	Now Lemma \ref{ringdecomp} implies that our graph $G$ has a ring-decomposition of width at most $w$ and length at least $n^{2^p}$, where $p = \lceil \log_2 w\rceil$. Note that by absorbing any terms $G_i$ where $i>n^{2^p}$ into $G_0$, we can ensure that our ring-decomposition has length exactly $n^{2^p}$. By our previous assertion, either $G$ has a connected ring-decomposition of width at most $w$ and length at least $n^{2^{p-1}}$, or $G$ has a ring-decomposition of width at most $w/2$ and length at least $n^{2^{p-1}}$. In the first case, the theorem holds. In the second case, we apply the previous assertion again to either find a connected ring-decomposition of sufficient length or to find a ring-decomposition of yet smaller width. By our choice of $p$, we guarantee that this process will terminate, yielding a connected ring-decomposition (possibly having width $1$) of length at least $n$. \hfill $\Box$ \\
	
	In the next section, we will see how to use connected ring-decompositions to find the immersions we desire.

	\section{Finding double cycles}
	
	Suppose that $A$ and $B$ are disjoint sets of vertices of a graph $G$. Let $Q_1 = x_1, \ldots, x_s$ and $Q_2 =y_1,\ldots, y_t$ be edge-disjoint $AB$-paths, where $\{x_1,y_1\}\subseteq A$ and $\{x_s, y_t\}\subseteq B$. Suppose $Q_1$ and $Q_2$ intersect in two internal vertices $u = x_{i_1} = y_{j_1}$ and $v = x_{i_2} = y_{j_2}$. If $(i_1 - i_2)(j_1-j_2) > 0$, then we say that $Q_1$ and $Q_2$ are \textit{aligned with respect to $u$ and $v$}. If $Q_1$ and $Q_2$ are aligned with respect to every pair of internal vertices in which they intersect, we say that the two paths are \textit{aligned}. A collection of paths which are pairwise aligned is also called aligned. 
	
	\begin{lemma}\label{aligned}
		For any collection $\mathcal{Q}$ of $k$ edge-disjoint $AB$-paths, the union of the members of $\mathcal{Q}$ contains a collection of $k$ edge-disjoint $AB$-paths which is aligned.
	\end{lemma}
	
	\noindent\textit{Proof:} Let $\mathcal{Q} = \{Q_1,\ldots, Q_k\}$ be a collection of $k$ edge-disjoint $AB$-paths. Among all such collections contained in $\bigcup\limits_{i=1}^k E(Q_i)$, let $\mathcal{Q}' =\{Q_1', \ldots, Q_k'\}$ denote a collection such that $\sum\limits_{i=1}^k |E(Q_i')|$ is minimized. We claim that $\mathcal{Q}'$ is aligned. 
	
	Suppose that $\mathcal{Q}'$ is not aligned. Then it contains some pair, which we may assume are $Q_1'$ and $Q_2'$, that are unaligned with respect to some pair of internal vertices $u$ and $v$. We may write $Q_1' = x_1, \ldots, u, \ldots, v\ldots, x_s$ and $Q_2' = y_1,\ldots, v, \ldots, u,\ldots, y_t$. Then by defining the paths $Q_1''= Q_1'[x_1,u]\cup Q_2'[u,y_t]$ and $Q_2''= Q_2'[y_1,v]\cup Q_1'[v,x_s]$, we find that the collection $\{Q_1'', Q_2'',\ldots, Q_k'\}$ forms a collection of $k$ edge-disjoint $AB$-paths with fewer edges than $\mathcal{Q}'$, since the edges of the subpaths $Q_1'[u,v]$ and $Q_2'[v,u]$ are omitted. This contradiction implies that $\mathcal{Q}'$ is aligned, as desired. \hfill $\Box$ \\
	
	For vertices $x,y$ in a graph $G$, we will call a set of edge-disjoint $xy$-paths $r$-{\it{compatible}} if some two of the paths are aligned and intersect in at least $r$ vertices (including $x$ and $y$). Before we begin to prove the main result of this paper, we make the following observation. 
	
	\begin{obs} \label{compatible}
		If a graph $G$ contains an $r$-compatible set of 4 $xy$-paths for some $x,y\in V(G)$, then $G$ contains an immersion of $C_{2,r}$.
	\end{obs} 
	
	\noindent\textit{Proof}: Given a set of $r$-compatible $xy$-paths $\{Q_1, Q_2, Q_3, Q_4\}$ in $G$, we may assume that $Q_1$ and $Q_2$ intersect in the vertices $\{x=v_1, v_2, \ldots, v_r = y\}$. Then, since $Q_1$ and $Q_2$ are aligned, we may lift $Q_j[v_i, v_{i+1}]$ for $i\in[r-1]$ and $j\in[2]$ to yield a $P_{2,r}$ immersion in $G$. Lifting the entirety of $Q_3$ and $Q_4$ completes a $C_{2,r}$ immersion in $G$. \hfill $\Box$ \\
	
	In trying to prove our main result, we often find two $xy$-paths which suffice to create a $P_{2,r}$ immersion of a graph $G$, and we require two additional paths between $x$ and $y$ to complete the $C_{2,r}$ immersion we desire. The following lemma provides these paths. For a subgraph $H$ and a vertex $x$ of a graph $G$, let $H\cup x$ denote the subgraph of $G$ consisting of $H$, $x$, and all edges between $x$ and $V(H)$.
	
	\begin{lemma} \label{Hplus}
		Let $G$ be a graph containing edge-disjoint connected subgraphs $G_1, G_2, G_3,$ and $G_4$ such that $H := G_1\cup G_2$ and $K := G_3\cup G_4$ are disjoint and $G_1\cap G_2$ and $G_3\cap G_4$ are nonempty. Suppose further that $G$ contains the following six edge-disjoint paths:
		\begin{itemize}
			\item $L_i^a$ between $a_i\in G_1 - G_2$ and $b_i \in G_2 - G_1$ for $i=1,2$,
			\item $L_i^c$ between $c_i\in G_3 - G_4$ and $d_i\in G_4 - G_3$ for $i=1,2$, and
			\item $L_3$ and $L_4$ such that either:
			\begin{enumerate}[a.] \item $L_3$ and $L_4$ are between $G_1\cap G_2$ and $G_3\cap G_4$, or
				\item $L_3$ is between $G_1$ and $G_4$ and $L_4$ is between $G_2$ and $G_3$ and $L_3\cup L_4$ is edge-disjoint from $H\cup K$.
			\end{enumerate}
		\end{itemize}
		
		Lastly, suppose $G$ contains vertices $x$ and $y$ (not contained in any $G_i$ or any of the six paths above) of degree two such that the neighbors of $x$ are $b_1$ and $b_2$ and the neighbors of $y$ are $c_1$ and $c_2$. Then $G$ contains two edge-disjoint $xy$-paths.
	\end{lemma}
	
	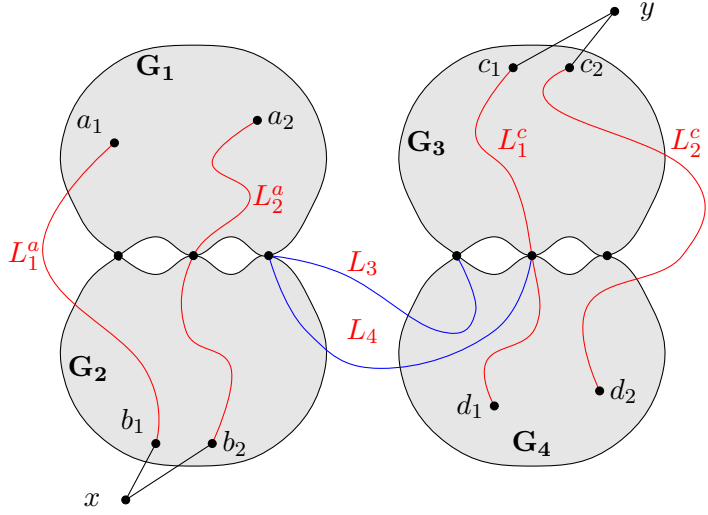
\begin{figure}[ht!]
		\begin{center}
			\begin{tikzpicture}[scale=0.5]
				
				\filldraw[black!10,draw=black] plot [smooth cycle,tension=.75] coordinates {(0,0) (-1,-1) (-1.5,-3) (-.25,-5) (2,-5.6) (4.25,-5) (5.5,-3) (5,-1) (4,0) (3,-.5) (2,0) (1,-.5) (0,0)};

				\begin{scope}[yscale=-1] 
					\filldraw[black!10,draw=black] plot [smooth cycle,tension=.75] coordinates {(0,0) (-1,-1) (-1.5,-3) (-.25,-5) (2,-5.6) (4.25,-5) (5.5,-3) (5,-1) (4,0) (3,-.5) (2,0) (1,-.5) (0,0)};

				\end{scope}
				
				\begin{scope}[shift = {(9,0)}]
					
					\filldraw[black!10,draw=black] plot [smooth cycle,tension=.75] coordinates {(0,0) (-1,-1) (-1.5,-3) (-.25,-5) (2,-5.6) (4.25,-5) (5.5,-3) (5,-1) (4,0) (3,-.5) (2,0) (1,-.5) (0,0)};

					\begin{scope}[yscale=-1]
						\filldraw[black!10,draw=black] plot [smooth cycle,tension=.75] coordinates {(0,0) (-1,-1) (-1.5,-3) (-.25,-5) (2,-5.6) (4.25,-5) (5.5,-3) (5,-1) (4,0) (3,-.5) (2,0) (1,-.5) (0,0)};
					\end{scope}
				\end{scope}
				
				\coordinate (x) at (.2, -6.5); \coordinate (u1) at (1, -5); \coordinate (u2) at (2.5, -5); 
				\coordinate (y) at (13.2, 6.5); \coordinate (v1) at (12, 5); \coordinate (v2) at (10.5, 5); 
				\coordinate (a) at (-.1, 3); \coordinate (b) at (3.7, 3.6); \coordinate (c) at (10, -4); \coordinate (d) at (12.8, -3.6);
				
				\draw (u1) -- (x)--(u2); \draw (v1) -- (y)--(v2);
				\draw[red] plot [smooth,tension=.75] coordinates {(u1) (.7,-3) (-2,0) (a)}; 
				
				\draw[red] plot [smooth,tension=.75] coordinates {(u2) (3,-3) (1.8,-2) (2,0) 
					(3.5,1.5) (2.5,2.5) (b)};
				
				\draw [blue] plot [smooth,tension=.75] coordinates {(4,0) (5,-2) (7,-3) (10,-2) (11,0)}; 
				
				\draw[blue] plot [smooth,tension=.75] coordinates {(4,0) (6,-.5) (8.5,-2) (9.5,-1.5) (9,0)}; 
				
				\draw[red] plot [smooth,tension=.75] coordinates {(c) (9.8,-2.8) (11.1,-1.9) (11,0) (10.5,2) (9.5,3.3) (v2)}; 
				
				\draw[red] plot [smooth,tension=.75] coordinates {(d) (12.5,-1.3) (15,0) (15.2,2.2) (11.5,4) (v1)};

				\foreach \coord in {x,u1,u2,y,v1,v2,a,b,c,d}{
					\draw[fill=black] (\coord) circle (3pt);}
				
				\foreach \x in {0,2,4,9,11,13}{
					\draw[fill=black] (\x,0) circle (3pt);}
				
				\draw (x) node[left=.2cm] {$x$}; \draw (y) node[right=.2cm] {$y$};
				\draw[red] (-2.5, 0) node {$L_1^a$}; \draw[red] (4, 1.5) node {$L_2^a$};
				\draw[red] (6.5, -.2) node {$L_3$}; \draw[red] (6.5, -2) node {$L_4$};
				\draw[red] (10.5, 3) node {$L_1^c$}; \draw[red] (15.1, 3) node {$L_2^c$};
				\draw (u1) node[above left] {$b_1$}; \draw (u2) node[right] {$b_2$};
				\draw (a) node[above left] {$a_1$}; \draw (b) node[right] {$a_2$};
				
				\draw (v2) node[left] {$c_1$}; \draw (v1) node[right] {$c_2$};
				\draw (c) node[left] {$d_1$}; \draw (d) node[right] {$d_2$};
				
				\draw (1,5) node {$\mathbf{G_1}$}; \draw (-.8,-3) node {$\mathbf{G_2}$}; 
				\draw (8.2,3) node {$\mathbf{G_3}$}; \draw (11,-5) node {$\mathbf{G_4}$}; 
				
			\end{tikzpicture}
		\end{center}
		\caption{A visual aid to Lemma \protect\ref{Hplus}(a).}
	\end{figure}
	
	\noindent\textit{Proof:} We may assume that $G$ is the union of the four $G_i$ and the six paths described above, along with the four edges incident to $x$ and $y$; it is clear that $G$ is connected. Suppose that $G$ does not contain two edge-disjoint $xy$-paths. Then by Menger's theorem, there exists a cut-edge between $x$ and $y$ in $G$. In other words, there exists a partition $(A,B)$ of $V(G)$ such that $x\in A,\, y\in B$, and exactly one edge $e$ of $G$ has one endpoint in $A$ and the other in $B$. Note that this implies that the vertices of any connected subgraph $G'$ of $G$ are split between $A$ and $B$ if and only if $e\in E(G')$.  We examine where in $G$ the edge $e$ could exist; due to the symmetry of $G$, we need only consider whether $e$ is incident to $x$, $e\in G_1$, $e\in G_2$, or $e\in (L_1^a\cup L_2^a \cup L_1^c \cup L_2^c \cup L_3 \cup L_4)\bs (H \cup K)$. \\
	
	\noindent\textbf{Case (i):} $e$ is incident with $x$ or $e\in G_1$ \\
	
	In this case we find $\{x\}\cup V(G_2) \subseteq A$. If (a.) holds, we may assume that $e\notin L_4$; if (b.) holds, we know this is the case. In either case, we find that $(G_2 \cup L_4 \cup K\cup y)$ is a connected subgraph of $G$ not containing $e$, and hence we must have that $y\in A$, a contradiction. \\
	
	\noindent\textbf{Case (ii):} $e\in G_2$ \\ 
	
	We may assume that $L_1^a\cup xb_1$ avoids $e$. Then  $x\cup L_1^a \cup G_1$ is a connected subgraph of $G$ not containing $e$, and hence all of the vertices of this subgraph are contained in $A$. If (a.) holds, we may assume that $e\notin L_3$; if (b.) holds, we know this is the case. In either case, we find that $x\cup L_1^a \cup G_1 \cup L_3 \cup K \cup y$ is a connected subgraph of $G$ not containing $e$, and hence $y\in A$, a contradiction. \\
	
	\noindent\textbf{Case (iii):} $e\in (L_1^a\cup L_2^a \cup L_1^c \cup L_2^c \cup L_3 \cup L_4)\bs (H\cup K)$ \\ 
	
	In this case, we see immediately that $\{x\}\cup V(G_1) \cup V(G_2)\subseteq A$. But here we must have that $e$ is in at most one of $L_3$ and $L_4$, and hence we may assume (in either case (a.) or (b.)) that $e\notin L_3$. Then we find that $x\cup H \cup L_3 \cup K \cup y$ is a connected subgraph of $G$ not containing $e$, which implies $y\in A$, a contradiction. \\
	
	Since we have reached a contradiction in each case, we see that no such edge $e$ exists in $G$, and hence our lemma holds. \hfill $\Box$\\
	
	Another useful tool is the following lemma, which allows us to perform contractions while maintaining edge-connectivity. For a set $S\subseteq V(G)$, we define the contraction of $S$ as the graph $G\slash S$ obtained by identifying all of the vertices in $S$ into a new vertex $x$ and then deleting any resulting loops. Note that if $G$ contains $k$ edges between $S$ and some vertex $y\notin S$, then the contraction of $S$ will result in $k$ edges in parallel between $x$ and $y$. The following lemma follows from the fact that every edge-cut in $G\slash S$ is an edge-cut of $G$.
	
	\begin{lemma} \label{contraction}
		Let $G$ be $k$-edge-connected and let $S\subset V(G)$. Then $G\slash S$ is $k$-edge-connected. 
	\end{lemma}
	
	Next, we prove a rather specific lemma, the utility of which will become evident later. In what follows, it is important to clarify that in any lifting involving the edges $e_i$ or $f_i$, we will call the newly created edges $e_i$ or $f_i$, respectively. 
	
	\begin{lemma}\label{claim}
		
		Let $G$ be a 4-edge-connected graph with $\Delta(G) \leq 5$ and let $A\subseteq V(G)$ such that $\delta(A) = \{e_1, \ldots, e_w, f_1,\ldots, f_w\}$ is a matching. If $H:=G\bs E(G-A)$  contains $w$ vertex-disjoint paths $R_i$ ($i\in[w]$), where each $R_i$ joins $e_i$ and $f_i$, then $H$ contains an immersion $H'$ which contains $w$ vertex-disjoint paths $Q_i$ ($i\in[w]$) joining $e_i$ and $f_i$ and a cycle $C$ such that $E(C\cap Q_i) = \emptyset$ for all $i\in[w]$ but $V(C\cap Q_s) \neq\emptyset \neq V(C\cap Q_t)$ for some distinct  $s,t\in[w]$.
	\end{lemma}
	
	\noindent\textit{Proof:} Let $S = V(G)\bs V(A)$, and let $H^+$ denote the graph obtained from $G$ by contracting the set $S$ into a single vertex $z$. By Lemma \ref{contraction}, we have that $H^+$ is $4$-edge-connected. 
	
	Note that finding the desired immersion in $H$ is equivalent to finding in $H^+$ an immersion containing $w+1$ cycles $C, D_1, \ldots, D_w$ such that each $D_i$ contains $e_i$ and $f_i$, $V(D_i\cap D_j) = \{z\}$ for all distinct $i,j\in[w]$, and $E(C\cap D_i) = \emptyset$ for all $i\in[w]$ but $V(C\cap D_s) \neq\emptyset \neq V(C\cap D_t)$ for some distinct  $s,t\in[w]$. Since we are looking for an immersion of $H^+$, we may assume that our $H^+$ is immersional-minimal with respect to the following three properties:
	
	\begin{enumerate}
		\item Every vertex of $H^+$ other than $z$ has degree 4 or 5.
		\item $H^+$ is $4$-edge-connected.
		\item $H^+$ has $w$ cycles $D_1, D_2, \ldots, D_{w}$ such that each $D_i$ contains $e_i$ and $f_i$ and $V(D_i\cap D_j) = \{z\}$ for all distinct $i,j\in[w]$. Let $\mathcal{D} = D_1\cup D_2\cup \ldots \cup D_{w}.$ 
		
	\end{enumerate}
	
	Let us further assume that $H^+$ is chosen among all such minimal immersions so that $|V(\mathcal{D})|$ is minimized. We consider additional properties of such an immersion.
	
	It is clear from our choice of $H^+$ that each cycle $D_i$ is an induced cycle (possibly with parallel edges). \\
	
	\noindent \textbf{Claim 1:} Every vertex outside of $V(\mathcal{D})$ has degree 5, and every edge $e=uv$ outside of $E(\mathcal{D})$ has an end with degree 4. \\
	
	\noindent \textit{Proof:} Note that any vertex outside of $V(\mathcal{D})$ having degree 4 could be eliminated via Mader's Theorem. Suppose that both $u,v$ have degree 5. Then the graph $H^+\bs e$ is an immersion of $H^+$ which clearly satisfies properties (1) and (3). This means that $H^+\bs e$ must admit a $3$-edge-cut; that is, $H^+$ admits a $4$-edge-cut containing $e$. This cut separates $H^+$ into two components, say $Z$ and $Y$, such that $z\in V(Z)$. Consider the graph $H'$ obtained from $H^+$ by contracting $V(Y)$ to a single vertex. It is clear that $H'$ satisfies property (1). Furthermore, Lemma \ref{contraction} implies that $H'$ satisfies (2). Lastly, since at most one of the $D_i$ meets the edge-cut containing $e$, we have that $H'$ also satisfies property (3).
	
	We now claim that $H'$ is a proper immersion of $H^+$, violating the minimality of our choice of $H^+$. To see this, choose a vertex $y\in V(Y)$. Then there are $4$ edge-disjoint paths from $y$ to $Z$, and these must travel through the $4$-edge-cut. After lifting the entirety of these $4$ paths, the resulting graph contains a subgraph isomorphic to $H'$, implying that $H'$ is an immersion of $H^+$. Because both ends of $e$ have degree 5, there must be at least one edge with both ends in $Y$, implying that $H'$ is a proper immersion of $H^+$. \hfill $\Box$ \\
	
	Recall the notion of a bridge defined in Section 2. Claim 1 implies that every trivial bridge of $\mathcal{D}$ is incident with a vertex of degree 4. It also implies that $H^+- \mathcal{D}$ is edgeless, which we can use to see the following about nontrivial $\mathcal{D}$-bridges in $H^+$: \\
	
	\noindent \textbf{Claim 2:} Every nontrivial bridge $B$ of $\mathcal{D}$ is a $K_{1,t}$ plus parallel edges, such that its center vertex has degree 5 and its feet have degree 4 in $H^+$. Furthermore, for any fixed $i\in[w]$, $B$ has either one foot on $D_i$, or exactly two feet on $D_i$, which are either adjacent or of distance 2 and separated by a vertex of degree 5. \\
	
	\noindent \textit{Proof:} Claim 1 immediately yields the fact that $B$ is a $K_{1,t}$ plus parallel edges, such that its center vertex has degree 5 and its feet have degree 4 in $H^+$. Fix an $i\in[w]$. Then $B$ has at most three feet on $D_i$ by the minimality of $\mathcal{D}$; otherwise, $D_i$ could be rerouted through the nucleus of $B$. In fact, $B$ does not have three feet on $D_i$, since these would have to be consecutively adjacent. In this case, the middle foot could be removed via Mader's theorem and replaced by the nucleus of the bridge, yielding an immersion satisfying properties (1)-(3) and violating the immersion-minimality of our choice of $H^+$. This means that $B$ has either at most one foot on $D_i$, or exactly two feet, which are either adjacent or of distance 2 and separated by a vertex of degree 5. \hfill $\Box$ \\
	
	Let $K_1, K_2, \ldots$ be components of $H^+ \bs \mathcal{D}$. Note that each $K_i$ is a union of bridges of $\mathcal{D}$, and that distinct $K_i, K_j$ are vertex disjoint. Let us fix one of these components and call it $K$. Let $F$ be the set of cut-edges of $K$. Then components of $K\bs F$ are maximal $2$-edge-connected subgraphs of $K$; let us call these components \textit{blobs} of $K$. Let $B_0$ be a leaf blob, a blob that is incident with only one edge of $F$. 
	
	First we see that $B_0$ is not a single vertex. Suppose that $xy$ is the cut-edge of $K$ with $x\in B_0$. Since $|\delta(x)| \geq 4$ and $\delta(x)$ contains at most two edges of $\mathcal{D}$, we must have that at least one edge of $\delta(x)$ belongs to $B_0$. But this implies that $B_0$ contains at least one vertex of $\mathcal{D}$, by the fact that $H^+ - \mathcal{D}$ is edgeless. We may assume that $B_0$ contains a vertex from $D_1$. If $B_0$ also contains a vertex from $\mathcal{D} - D_1$, then $B_0$ contains the desired cycle. 
	
	Suppose that this is not the case. Then $B_0\cap \mathcal{D} = B_0 \cap D_1 \neq \emptyset$, which further implies that $V(B_0) \subseteq V(D_1)$. Otherwise, $B_0$ contains a vertex $x$ of $H^+\bs\mathcal{D}$, which implies that $x$ has degree 5. The choice of $B_0$ requires that at most one edge in $\delta(x)$ is a cut-edge of $K$, which implies that at least four edges of $\delta(x)$ are in $E(B_0)$. However, by Claim 2, $x$ has exactly two neighbors on $D_1$ (a unique neighbor would have degree greater than 5). Then either we see at least two edges of $\delta(x)$ incident with vertices of $\mathcal{D}\bs D_1$, a contradiction, or we have that both of these neighbors (which must have degree 4) are connected to $x$ via a parallel pair of edges. In the latter case, contracting $x$ and one of its neighbors into a single vertex yields an immersion (via Lemma \ref{contraction}) that violates our choice of $H^+$.
	
	The fact that $V(B_0) \subseteq V(D_1)$ implies that $B_0$ consists of two edges in parallel to an edge of $D_1$. To see this, note that our previous argument and the fact that the $D_i$ are induced cycles with parallel edges together imply that all edges of $B_0$ are parallel to edges of $D_1$. Since every vertex of $H^+$ other than $z$ has degree less than 6, we must have that $B_0$ is as described. However, this yields an edge-cut $\delta(V(B_0))$ of size 3 in $H^+$, a contradiction. Therefore $B_0$ must contain a vertex from $\mathcal{D}\bs D_1$, and hence contains the desired cycle. \hfill $\Box$  \\
	
	The preceding lemmas form a powerful toolbox for proving the following lemma, which is the most important step in proving our main result.
	
	\begin{lemma}\label{doublecycle}
		For any integers $n\geq 4$ and $w\geq 1$, there exists a number $\ell = f_{\ref{doublecycle}}(n,w)$ such that any $4$-edge-connected graph $G$ admitting a connected ring-decomposition of width $w$ and length at least $\ell$ contains $C_{2,n}$ as an immersion.
	\end{lemma}
	
	\noindent\textit{Proof:} We introduce the following intermediate numbers which appear in this proof: 
	\begin{align*}
		n_1 &= \binom{w+1}{2}(n-1) +5 \\
		n_2 &= \binom{w}{2}(n-1) + 6 \\
		\ell &= f_{\ref{doublecycle}}(n,w) = n_1n_2
	\end{align*} 
	
	Suppose that $G$ is a $4$-edge-connected graph admitting a connected ring-decomposition $(G_0, G_1, \ldots, G_\ell)$ of width $w$, and let $G' = G_1 \cup \ldots \cup G_\ell$. The proof is divided into cases depending upon the number of edge-disjoint $W_1W_{\ell+1}$-paths in $G'$, where $W_i = V(G_{i-1}\cap G_i)$. \\
	
	\noindent\textbf{Case 1:} $G'$ contains $w+1$ edge-disjoint $W_1W_{\ell+1}$-paths $\{Q_1, \ldots, Q_{w+1}\}$. \\
	
	To settle this case, we need only that $\ell \geq n_1$. Let $X = W_3$ and $Y = W_{\ell-1}$. By Lemma \ref{aligned}, we may assume that the collection $\{Q_1, \ldots, Q_{w+1}\}$ is aligned. Suppose that each path $Q_i$ travels from $a_i\in W_1$ to $d_i\in W_{\ell+1}$. Let $b_i$ denote the last vertex of $Q_i$ in $X$ as one travels from $a_i$ to $d_i$, and let $c_i$ denote the first vertex of $Q_i$ in $Y$ as one travels from $b_i$ to $d_i$. Then the subpaths $Q_i^a = Q_i[a_i,b_i]$, $Q_i^b = Q_i[b_i,c_i]$, and $Q_i^c = Q_i[c_i,d_i]$ partition the edges of $Q_i$. Furthermore, by the choice of $b_i$ and $c_i$, we see that $Q_i^b$ is contained in $G_3\cup\ldots\cup G_{\ell-2}$. In fact, we see that the collection $\{Q_1^b,\ldots, Q_{w+1}^b\}$ forms an aligned set of $w+1$ edge-disjoint $XY$-paths.
	
	Because there are $w+1$ paths intersecting a gate of size $w$, the Pigeonhole Principle implies that in each $W_j$ ($3\leq j \leq \ell-1$), there must be some vertex belonging to at least two of the $Q_i^b$. That is, from the $\binom{w+1}{2}$ pairs of the paths $Q_i^b$, each $W_j$ can be assigned a number representing a pair of paths which intersect in one of its vertices. For any sequence of gates of length $\binom{w+1}{2}(n-1)+1$, then, there must be some pair of paths, say $Q_1^b, Q_2^b$, which intersect at least $n$ times. Furthermore, for any sequence of gates of length $n_1 = \binom{w+1}{2}(n-1)+5$, we may assume that these intersections occur in the stretch between $X$ and $Y$. Since $\ell \geq n_1$, we know that such a sequence exists in $G'$. Let $v_1, v_2, \ldots, v_n$ be the vertices of intersection; these will form the terminals of our immersed $C_{2,n}$. Note that $Q_1[b_1, v_1]$ and $Q_2[b_2, v_1]$ provide two edge-disjoint paths from $v_1$ to $X$, while $Q_1[v_n, c_1]$ and $Q_2[v_n, c_2]$ provide two edge-disjoint paths from $v_n$ to $Y$.\\

	\noindent\underline{Subcase 1(a):} $w\geq 2$ \\
	
	In this case $G'$ contains a set of 3 edge-disjoint $W_1W_{\ell+1}$-paths $\{Q_1, Q_2, Q_3\}$ and, consequently, a set of 3 aligned edge-disjoint paths $\{Q_1^b, Q_2^b, Q_3^b\}$ between $X$ and $Y$ such that $Q_1^b$ and $Q_2^b$ intersect in $n$ vertices. We claim that these can be extended to form a set of $4$ edge-disjoint $v_1v_n$-paths which are $n$-compatible. To do so, we need only provide two $v_1v_n$-paths which are edge-disjoint from each other and also from $Q_1^b[v_1,v_n]\cup Q_2^b[v_1,v_n]$; we will utilize Lemma \ref{Hplus}.
	
	Consider the number of edge-disjoint $W_1W_{\ell+1}$-paths in $G'$. If this number is at least four, then $G'$ contains a $W_1W_{\ell+1}$ path $Q_4$ which is edge-disjoint from the $Q_i^b$ $(i\in[3]$). If this number is less than four, then we must have that $w=2$ and $G'$ contains exactly three edge-disjoint $XY$-paths, so that $G'$ contains an edge-cut of size three between $W_1$ and $W_{\ell+1}$. But then the 4-edge-connectivity of $G$ and Menger's theorem together imply that $G_0$ contains a $W_1W_{\ell+1}$-path.
	
	To apply Lemma \ref{Hplus}, we first lift each of the paths $Q_1[b_1, v_1]$, $Q_2[b_2, v_1]$, $Q_1[v_n, c_1]$ and $Q_2[v_n, c_2]$, and then make the following substitutions:
	
	\begin{itemize}
		\item $G_1 = G_1$ and $G_2 = G_2$,
		\item $G_3 = G_{\ell-1}$ and $G_4 = G_\ell$,
		\item $L_i^a$ = $Q_i^a$ and $L_i^c = Q_i^c$ for $i=1,2$,
		\item if $w>2$, let $L_i= Q_i'$ (for $i=3,4$) where $Q_i'$ is a subpath of $Q_i$ between $W_2$ and $W_l$; apply Lemma \ref{Hplus}(a).
		\item if $w=2$, let $L_3 = Q_3^b$ and $L_4$ be a $W_1W_{l+1}$-path in $G_0$; apply Lemma \ref{Hplus}(b).
		\item $x = v_1$ and $y = v_n$.
	\end{itemize}
	
	Lemma \ref{Hplus} implies that $G$ contains two edge-disjoint $v_1v_n$-paths which are further edge-disjoint from both $Q_1[v_1,v_n]$ and $Q_2[v_1,v_n]$. Thus we obtain a set of 4 edge-disjoint $v_1v_n$ paths which are $n$-compatible, and hence $G$ contains a $C_{2,n}$ immersion. \hfill $\Box$ \\
	
	\noindent\underline{Subcase 1(b):} $w = 1$ \\ 
	
	In this case, the graphs $G_i$ for $i\in[\ell]$ are joined in sequence by 1-sums, with one rail running from $W_1 = \{x\}$ to $W_{\ell+1} = \{y\}$. Note that any $xy$-paths in $G'$ are funneled through the gates and hence intersect in the unique vertex of $W_i$ for  $i\in \{2,\ldots, \ell\}$. Then:
	\begin{itemize}
		\item If $G'$ contains at least four edge-disjoint $xy$-paths, then this set of paths forms a set of four $n$-compatible $xy$-paths.
		\item If $G'$ contains exactly three edge-disjoint $xy$-paths, then $G'$ must contain a $3$-edge-cut between $x$ and $y$ by Menger's theorem. This edge-cut implies that $x$ and $y$ occupy the same component of $G_0$, guaranteeing an $xy$-path in $G_0$; this path along with the three in $G'$ forms a set of four $n$-compatible $xy$-paths.
		\item If $G'$ contains exactly two edge-disjoint $xy$-paths, then $G'$ must contain a $2$-edge-cut between $x$ and $y$ by Menger's theorem. Again, the existence of this cut implies that $G_0$ contains two edge-disjoint $xy$-paths; these paths along with the two in $G'$ forms a set of four $n$-compatible $xy$-paths.
	\end{itemize}
	
	\noindent\textbf{Case 2:} $G'$ contains exactly $w$ edge-disjoint $W_1W_{\ell+1}$-paths $\{Q_1, \ldots, Q_w\}$. \\ 
	
	Suppose that $w=1$. Then by Menger's theorem, $G'$ contains a unique cut-edge $x'y'$ between $W_1$ and $W_{\ell+1}$, since two such edges would yield a 2-edge-cut of $G$. Since $\ell > 2n$, we can assume (by symmetry) that $\{x'\} = W_{\ell'}$ for some $\ell'\geq n$ and that there exist at least two edge-disjoint $W_1W_{\ell'}$-paths in $G'$. Therefore $G$ contains a $C_{2,n}$ immersion by Subcase 1(b). 
	
	We proceed assuming that $w\geq 2$. Here, Menger's theorem implies that $G'$ contains a $w$-edge-cut $F$ between $W_1$ and $W_{\ell+1}$. However, our ring-decomposition guarantees the existence of $w$ vertex-disjoint rails $\{R_1, \ldots, R_w\}$ from $W_1$ to $W_{\ell+1}$, and so we may take $Q_i=R_i$ for each $i\in[w]$. This implies that $F$ is comprised of exactly one edge from each rail. 
	
	We would now like to produce a long sequence of disjoint edge-cuts in $G'$. To this end, we may assume that the length of our ring-decomposition is exactly $\ell$ by Remark \ref{absorb}. Note that if any sequence of $n_1$ consecutive $G_i$ contains no $w$-edge-cut, then that sequence contains $w+1$ edge-disjoint paths from its first gate to its last, which implies that $G$ contains a $C_{2,n}$ immersion by Case 1. Hence we may assume that any consecutive sequence of $n_1$ of the $G_i$ contains a $w$-edge-cut. For $i\in [n_2]$, let $G_i' = G_{n_1(i-1)+1} \cup \ldots \cup G_{in_1}$. Then $\{G_0, G_1', \ldots, G_{n_2}'\}$ forms a connected ring-decomposition of width $w$ and length $n_2$ such that a $w$-edge-cut $F_i$ separates the first and last gates of $G_i'$. Let $W_1' := W_1, W_{n_2+1}' := W_{n_1n_2+1}$, and $W_i' := G_i'\cap G_{i-1}'$ for $i\in \{2,\ldots, n_2\}$. Note that a given $G_i'$ may contain more than one $w$-edge-cut; we choose one arbitrarily. This yields a set of disjoint $w$-edge-cuts $\{F_1, F_2, \ldots, F_{n_2}\}$. Then, for $i\in[n_2-1],$ let $H_i$ denote the component of $G'\bs (F_i\cup F_{i+1})$ between the two edge cuts.
	
	We now have a sequence of disjoint edge-cuts along our $w$ rails $\{R_1,\ldots, R_w\}$ within $G'$, with $F_1, \ldots, F_{n_2}$ forming edge-disjoint subgraphs $H_i$ ($i\in[n_2-1]$) separated by these edge-cuts. It is this situation for which Lemma \ref{claim} prepares us. We may apply Lemma \ref{claim} to the subgraph $F_i \cup H_i \cup F_{i+1}$ to find that in each $H_i$, we may assume (by taking immersions) that we have a distinguished set of $w$ vertex-disjoint paths $\{Q^i_1, \ldots, Q^i_w\}$ and a distinguished cycle $C_i$ between some pair of the $Q^i_j$. For $j\in[w]$, let $R_j'$ denote the path obtained from $R_j$ by replacing the edges of $R_j$ within each $H_i$ by $Q_i^j$. Then $\{R_1', \ldots, R_w'\}$ is a set of $w$ vertex-disjoint paths from $W_1'$ to $W_{n_2+1}'$ such that in each $H_i$, there is a cycle between some two of the $R_j'$. Because $n_2 > \binom{w}{2}(n-1)+2$, there must be at least $\binom{w}{2}(n-1)+1$ of the $C_i$ and hence some pair of the $R_j'$, say $R_1'$ and $R_2'$, must be joined by some $C_i$ at least $n$ times. Let $C_1', \ldots, C_n'$ denote $n$ of the $C_i$ which join $R_1'$ and $R_2'$.  Then let $v_i = R_1' \cap C_i'$ and $u_i = R_2' \cap C_i'$ for $i\in[n]$; by lifting the two paths of $C_i'$ between $v_i$ and $u_i$, we can consider each $C_i'$ as a pair of parallel edges. We arbitrarily assign labels to the edges of each $C_i'$, say $C_i' = \{e_i, f_i\}$. Note that because $n_2 =\left(\binom{w}{2}(n-1)+2\right)+4$, we may assume that all of the $C_i$ (and hence $v_i$ and $u_i$) are contained in $H_3 \cup \ldots \cup H_{n_2-3}$.
	
	We will take the set of vertices $\{v_1,\ldots v_n\}$ to be the terminals of our desired immersion. Note that $R_1'$ provides one path joining each pair $v_i, v_{i+1}$ for $i\in[n-1]$. To find a $P_{2,n}$ immersion, we can lift the path $f_i \cup R_2'[u_i, u_{i+1}] \cup e_{i+1}$ for each $i\in[n-1]$. 
	
	The process of completing the $C_{2,n}$ immersion depends on the value of $w$. Let $w_i^j$ denote the unique vertex of $W_i'\cap R_j'$. 
	
	If $w=2$, then Menger's theorem implies that $G\bs (H_3 \cup \ldots \cup H_{n_2-3})$ contains two edge-disjoint paths between $F_3$ and $F_{n_2-3}$, which can be used to complete our $C_{2,n}$ immersion.
	
	If $w\geq 3$, we use Lemma \ref{Hplus} to complete our $C_{2,n}$ immersion. In the case that $w=3$, we find that $G'$ contains many 3-edge-cuts between $W_1'$ and $W_{n_2+1}'$, and the 4-edge-connectivity of $G$ and Menger's theorem together imply that $G_0$ contains at least one $W_1'W_{n_2+1}'$-path. 
	
	Now, by lifting the four paths $R_1'[v_1,w_3^1]$, $e_1\cup R_2'[u_1,w_3^2]$, $f_n \cup R_2'[u_n,w_{n_2-1}^2]$, and $R_1'[v_n,w_{n_2-1}^1]$, we can apply Lemma \ref{Hplus} where:
	\begin{itemize}
		\item $G_1 = G_1'$ and $G_2 = G_2'$;
		\item $G_3 = G_{n_2-1}'$ and $G_4 = G_{n_2}'$;
		\item $L_j^a = R_j'[w_1^j, w_3^j]$ and $L_j^c = R_j'[w_{n_2-1}^j, w_{n_2+1}^j]$ for $j=1,2$; 
		\item $x = v_1$ and $y = v_n$;
		
		\item if $w=3$, we apply Lemma \ref{Hplus}(b) where:
		\begin{itemize}[$\circ$]
			\item $L_3$ is a $W_1'W_{n_2-1}'$-path in $G_0$; and
			\item $L_4 = R_3'[w_3^3, w_{n_2-1}^3]$.
		\end{itemize}
		
		\item if $w\geq 4$, we apply Lemma \ref{Hplus}(a) where:
		\begin{itemize}[$\circ$]
			\item $L_3 = R_3'[w_2^3,w_{n_2}^3]$; and
			\item $L_4 = R_4'[w_2^4,w_{n_2}^4]$.
		\end{itemize}
		
	\end{itemize}
	
	Lemma \ref{Hplus} implies the existence of two $v_1v_n$-paths edge-disjoint from our current $P_{2,n}$ immersion; lifting each of these paths completes a $C_{2,n}$ immersion in $G$, as desired. \hfill $\Box$ \\
	
	Finally, we prove the main result of this paper. \\

	\noindent\textbf{Theorem \ref{main}} 
	\textit{For every integer $r\geq 4$, there exists a number $f_{\ref{main}}(r)$ such that every $4$-edge-connected graph of order at least $f_{\ref{main}}(r)$ contains $C_{2,r}$ as an immersion.} \\

	\noindent\textit{Proof:} For readability, we introduce the following intermediate numbers:
	\begin{align*}
		n_1 &= f_{\ref{chud}}(r) - 1 \\
		n_2 &= f_{\ref{doublecycle}}(r, n_1) \\
		f_{\ref{main}}(r) &= f_{\ref{connected}}(n_2, n_1,5)
	\end{align*} 
	
	Suppose that $G$ is a 4-edge-connected graph and that $|G| \geq f_{\ref{main}}(r)$. Recall that Observation \ref{obs1} allows us to assume that every vertex of $G$ has degree 4 or 5, and hence that $\Delta(G) \leq 5$. Corollary \ref{chud} states that $G$ contains a $C_{2,r}$ immersion if $tw(G) \geq f_{\ref{chud}}(r)$, so we may assume that $tw(G) \leq n_1$. This means that $G$ has bounded degree and bounded tree-width, so we turn to Theorem \ref{connected}. Given the order of our graph, Theorem \ref{connected} asserts that $G$ contains a connected ring-decomposition of width at most $n_1$ and length at least $n_2$. Finally, Lemma \ref{doublecycle} tells us that $G$ contains a $C_{2,r}$ immersion, finishing the proof. \hfill $\Box$.

	\section{Linear edge-connectivity forces \texorpdfstring{$C_{t,r}$}{Ctr} immersions}
	
	We would now like to discuss our second main result. The following is a result of B\"{o}hme, Kawarabayashi, Maharry, and Mohar \cite{bohme}.
	
	\begin{thm}
		For any integers $a, s,$ and $r$, there exists a constant $N(s, r, a)$ such that every $(3a + 2)$-connected graph of minimum degree at least $\frac{31}{2}(a+1)-3$ and with at least $N(s, r, a)$ vertices either contains $K_{a,sr}$ as a topological minor or a minor isomorphic to $s$ disjoint copies of $K_{a,r}$.
	\end{thm}
	
	By setting $a=t$ and $s=1$, we can state the following consequence of this result, which was mentioned in the introduction.\\
	
	\noindent\textbf{Theorem \ref{bigktk}} 
	\textit{A sufficiently large ($16t+13$)-connected graph contains a $K_{t,r}$ minor.} \\

	The desire to find an analogue of this result for edge-connectivity and unavoidable immersions motivated our next result.\\

	\noindent\textbf{Theorem \ref{ctr}} 
	\textit{For all positive integers $t,r$, there exists a number $f_{\ref{ctr}}(t,r)$ such that every $(6t-4)$-edge-connected graph of order at least $f_{\ref{ctr}}(t,r)$ contains $C_{t,r}$ as an immersion.} \\

	In proving this result, we will call upon the well-known fact that the lone unavoidable immersion of 2-edge-connected graphs is a long cycle, $C_r$. To our knowledge, however, no formal proof of this fact exists in the literature. We will provide one shortly. In doing so, we cite the following result, which appears in Diestel's \cite{diestel} text.
	
	\begin{prop}\label{2ctop}
		For every integer $r\geq 1$, there exists a number $f_{\ref{2ctop}}(r)$ such that every 2-connected graph of order at least $f_{\ref{2ctop}}(r)$ contains $C_r$ or $K_{2,r}$ as a topological minor.
	\end{prop}
	
	As mentioned, we state and prove the following for use later.
	
	\begin{prop} \label{2ecimm}
		For every integer $r\geq 3$, there exists a number $f_{\ref{2ecimm}}(r)$ such that every 2-edge-connected graph of order at least $f_{\ref{2ecimm}}(r)$ contains $C_r$ as an immersion.
	\end{prop}
	
	\noindent\textit{Proof:} Suppose that $G$ is a 2-edge-connected graph of order at least $f_{\ref{2ecimm}}(r) := f_{\ref{2ctop}}(r)$. Recall from Observation \ref{obs1} that we may assume that every vertex of $G$ has degree two or three. Note that a 2-edge-connected graph satisfying these degree conditions must be 2-vertex-connected. Therefore we may use Proposition \ref{2ctop} to see that $G$ contains either $C_r$ or $K_{2,r}$ as a topological minor. When $r=3$, it is straightforward to verify that $K_{2,3}$ immerses $C_3$. If $r>3$, we find $\Delta(G) \leq 3 < r$, so that must $G$ contain $C_r$ as a topological minor (and hence as an immersion). \hfill $\Box$ \\
	
	We also use the following result of Kundu \cite{kundu}, which itself is a corollary of a well-known result attributed to both Tutte \cite{tutte} and Nash-Williams \cite{nash}.
	
	\begin{thm} \label{kundu}
		A $2s$-edge-connected graph contains $s$ edge-disjoint spanning trees.
	\end{thm}
	
	The utility of this result is clear. By specifying the edge-connectivity of a graph $G$, we can guarantee the existence of as many edge-disjoint spanning trees of $G$ as our argument necessitates. Before we do so, however, we will establish a few tools. 
	
	Let us consider the following three graphs: $K_{1,t}, P_t,$ and $\text{comb}_t$. While the first two are ubiquitous, the third may need defining. To this end, we note that comb$_t$ is the graph obtained from $P_{t+2}$ by appending an edge to each internal vertex of the path.
	
	Our goal is to find structure within the spanning trees guaranteed by Theorem \ref{kundu}. Suppose that we have some specified set of vertices in a tree and that we would like to find some subtree containing many of our specified vertices. The following lemma allows us to do so, with even more specificity. Note that this is the finite analogue of Lemma 8.2.2 in Diestel's text \cite{diestel}, and we omit its straightforward proof.
	
	\begin{lemma}\label{comb}
		For every integer $t$, there exists a number $f_{\ref{comb}}(t)$ such that if $G$ is a connected graph with $X\subseteq V(G)$ and $|X|\geq f_{\ref{comb}}(t)$, then $G$ contains one of the following as a subgraph:
		\begin{enumerate}
			\item A path $P$ containing at least $t$ vertices in $X$;
			\item A subdivision of $K_{1,t}$ or comb$_t$ such that all of its leaves are in $X$.
		\end{enumerate}
	\end{lemma}
	
	The next lemma is the most important tool in proving Theorem \ref{ctr}. Before proving it, we wish to extend the idea of aligned paths to the intersections of different edge-disjoint subgraphs. To this end, we say that a path $P$ is aligned with the cycle $C$ if $P$ and $C$ are edge-disjoint and $P$ is aligned with $C\bs e$ for some edge $e$ of $C$. Next, let $K$ denote a comb$_t$ subgraph of a graph $G$ which has been obtained from a path $x_0,x_1,\ldots, x_{n+1}$ by adding pendant edges $\{x_1y_1,\ldots, x_ny_n\}$. Then $K$ is aligned with a cycle $C$ if $V(K\cap C) = \{x_0, y_1, \ldots, y_n, x_{n+1}\}$ and the (virtual) path $x_0y_1\ldots y_nx_{n+1}$ is aligned with $C$. Note that a subdivided comb will be called aligned with a cycle $C$ if the comb subgraph obtained from suppressing any vertices of degree 2 is aligned with $C$.
	
	\begin{lemma} \label{gauge}
		For any integer $n\geq 3$ there exists a number $m = f_{\ref{gauge}}(n)$ such that if a graph $G$ contains a $C_{p,m}$ subgraph $C$ and three edge-disjoint trees $T_a, T_b, T_c$, each of which span the vertices of $C$ and are edge-disjoint from $C$, then $C\cup T_a \cup T_b \cup T_c$ contains $C_{p+1,n}$ as an immersion. 
	\end{lemma}
	
	\noindent\textit{Proof:} We introduce the following intermediate numbers:
	\begin{align*}
		n_1 &= f_{\ref{comb}}(n^2) \\ 
		n_2 &= (n_1)^2 \\
		f_{\ref{gauge}}(n) &= f_{\ref{comb}}(n_2)
	\end{align*} 
	
	Let us write the vertex set of $C$ as $V(C) = X = \{x_1, x_2, \ldots, x_m\}$. By the choice of $m$, Lemma \ref{comb} implies that each $T_i$ contains either a path containing $n_2$ members of $X$, a subdivision of $K_{1,n_2}$ with all its leaves in $X$, or a subdivision of comb$_{n_2}$ with all its leaves in $X$. \\
	
	\noindent \textbf{Case 1:} Any one of $T_a, T_b,$ or $T_c$ contains a path containing at least $n^2$ members of $X$. \\
	
	We may assume it is $T_a$ that contains such a path, which we will call $P$. Let $Y = X \cap V(P)$. Since $|Y| \geq n^2$, we can choose a subset $Y' = \{y_1, \ldots, y_n\}$ of vertices in $Y$ such that if $P'$ is the path obtained from $P[y_1,y_n]$ by lifting the subpaths $P[y_i, y_{i+1}]$ for each $i\in [n-1]$, then $P'$ is aligned with $C$. Note that $C$ contains $p$ edge-disjoint paths between each $y_i$ and $y_{i+1}$; by lifting the entirety of these paths, we create an immersion of $C_{p+1, n} \bs e$ in $C\cup T_a$. Finally, $T_b$ contains a path between $y_1$ and $y_n$; by lifting this path, we find an immersion of $C_{p+1, n}$ in $C\cup T_a \cup T_b$, and the lemma holds. \\
	
	\noindent \textbf{Case 2:} None of $T_a, T_b,$ or $T_c$ contains a path containing $n^2$ members of $X$. \\
	
	In this case, we must have that $T_a$ contains a subdivision $K$ of either $K_{1,n_2}$ or comb$_{n_2}$, the leaves of which are in $X$. If $K$ is a subdivided comb, we can choose a subset $Y = \{y_1, \ldots, y_{n_1}\}$ of the leaves of $K$ such that the subcomb $K'$ of $K$ which has the members of $Y$ as its leaves is aligned with $C$, since $n_2 = (n_1)^2$. If $K$ is a subdivided star, we simply let $Y$ denote any subset of $n_1$ leaves of $K$ and select $K'$ as the substar of $K$ having $Y$ as its leaves. We now apply Lemma \ref{comb} to the set $Y$; since $n_1 = f_{\ref{comb}}(n^2)$, we must have that $T_b$ contains either a path containing $n^2$ vertices of $Y$ or a subdivision $J$ of either $K_{1,n^2}$ or comb$_{n^2}$ such that all the leaves of $J$ are in $Y$. Because $Y\subseteq X$, however, we have that $T_b$ does not contain a path containing $n^2$ members of $Y$. Therefore $T_b$ contains a subdivision $J$ of either $K_{1,n^2}$ or comb$_{n^2}$ such that all the leaves of $J$ are in $Y$. In the case that $J$ is comb, we can choose a subset $Z = \{z_1, \ldots, z_n\}$ of $Y$ such that the subcomb $J'$ of $J$ having $Z$ as its set of leaves is aligned with $C$. If $J$ is a subdivided star, we can let $Z = \{z_1, \ldots, z_n\}$ denote an arbitrary collection of $n$ vertices of $Y$ and let $J'$ denote the substar of $J$ having $Z$ as its set of leaves.
	
	The following is true no matter the structure of $K'$ and $J'$: $K'$ contains edge-disjoint paths $Q_i$ joining $z_i$ and $z_{i+1}$ for all odd $i$ in $[n-1]$, and $J'$ contains edge-disjoint paths $Q_i$ joining $z_i$ and $z_{i+1}$ for all even $i$ in $[n-1]$. We find, then, that $C \cup T_a \cup T_b$ contains an immersion of $C_{p+1, n} \bs e$. Finally, since $T_c$ contains a path between $z_1$ and $z_n$, we find an immersion of $C_{p+1, n}$ in $C\cup T_a \cup T_b\cup T_c$, and the lemma holds. \hfill $\Box$ \\
	
	With this lemma in hand, we can prove Theorem \ref{ctr}. For any function $f$, let $f^{(0)}(x) = x$. Then, for $n\geq 1$, let $f^{(n)}(x)$ denote the $n$-fold composition of $f$ (e.g. $f^{(2)}(x) = (f\circ f) (x)$). \\
	
	\noindent\textit{Proof of Theorem \ref{ctr}:} 
	Note that if $t=1$, this result is exactly Proposition \ref{2ecimm}, so we may take $f_{\ref{ctr}}(1,r) := f_{\ref{2ecimm}}(r)$. We continue under the assumption that $t\geq 2$. Suppose that $G$ is a $(6t-4)$-edge-connected graph with $|G| \geq  f_{\ref{ctr}}(t,r) := f_{\ref{main}}(f_{\ref{gauge}}^{(t-2)}(r))$. Theorem \ref{kundu} implies that $G$ contains $3t-2$ edge-disjoint spanning trees $T_1, T_2, \ldots, T_{3t-2}$. Note that $T_1\cup T_2 \cup T_3 \cup T_4$ is a $4$-edge-connected graph, and hence Theorem \ref{main} implies that it contains a double cycle immersion $C$ of length $f_{\ref{gauge}}^{(t-2)}(r)$. If $t=2$, then we are done. If $t\geq 3$,  the trees $T_5, T_6,$ and $T_7$ span the vertices of $C$, and by applying Lemma \ref{gauge}, we find that $G$ contains an immersion of $C_{3, f_{\ref{gauge}}^{(t-1)}(r)}$. We can continue in this fashion, applying Lemma \ref{gauge} a total of $t-2$ times in succession to find a $C_{t,r}$ immersion in $G$, as desired. \hfill $\Box$ \\
	
	It is also possible to prove a \textit{rooted} version of Theorem \ref{ctr}; that is, given a specified set of vertices $S$, we would like to choose a $C_{t,r}$ immersion so that all of the terminals of our immersion belong to $S$. We say that a set $S\subseteq V(G)$ is $k$-edge-connected if $S$ is never separated by the removal of less than $k$ edges from $G$. Furthermore, if some tree $T$ in $G$ contains all of the vertices of $S$, we call $T$ an \textit{$S$-Steiner-tree} (or simply $S$-tree) in $G$. Generalizing Theorem \ref{kundu}, Kriesell \cite{kriesell} conjectured that if $S$ is a $2k$-edge-connected set in $G$, then $G$ contains $k$ edge-disjoint $S$-trees. Lau \cite{lau} was the first to contribute towards this conjecture, proving that if $S$ is $24k$-edge-connected in $G$, then $G$ contains $k$ edge-disjoint $S$-trees. West and Wu \cite{west} later improved this bound to $6.5k$. To our knowledge, the best current bound is the following, due to DeVos, McDonald, and Pivotto \cite{DeVosbound}.
	
	\begin{thm}\label{DeVosbound}
		If $S$ is $5k+4$-edge-connected in $G$, then $G$ has $k$ edge-disjoint $S$-trees.
	\end{thm}
	
	This result allows us to prove the following theorem.
	
	\begin{thm}\label{rooted}
		For all positive integers $t,r$, there exists a number $f_{\ref{rooted}}(t,r)$ such that if $S\subseteq V(G)$ is $(15t-6)$-edge-connected and $|S| \geq f_{\ref{rooted}}(t,r)$, then $G$ contains a $C_{t,r}$ immersion having every terminal in $S$.
	\end{thm}
	
	Note that $15t-6 = 5(3t-2)+4$, so that $G$ contains $3t-2$ $S$-trees by Theorem \ref{DeVosbound}. We omit further details of the proof because of its similarity to the proof of Theorem \ref{ctr}. 
	
	\section{Immersions and Line Graphs}
	
	We conclude the paper with a few comments on immersions and line graphs. The following proposition is a natural one, and some mention of it is made in the introduction of \cite{chudnov}. For completeness, we state it formally and include a proof here.
	
	\begin{prop} \label{minor}
		Let $G$ and $H$ be connected graphs. If $H$ is immersed in $G$, then $L(H)$ is a minor of $L(G)$.
	\end{prop}
	
	\noindent\textit{Proof:} This result is most easily observed by considering the operations performed on $G$ to arrive at $H$ and considering their effects on $L(G)$. Recall that $H$ can be obtained from a subgraph of $G$ by repeatedly lifting pairs of edges. This means we may consider three operations -- deleting edges, deleting isolated vertices, and lifting pairs of edges -- and their effects on $L(H)$ and $L(G)$. First, deleting an edge of $G$ corresponds to deleting a vertex of $L(G)$, which is clearly an allowable operation in finding a minor. Next, we find that isolated vertices in $G$ simply do not appear in $L(G)$, and hence may be deleted without affecting $L(G)$. 
	
	Last, we consider lifting a pair of edges $\{uv,vw\}$. Let $\delta^*(v)$ denote the subset of $\delta(v)$ consisting of edges not incident to either $u$ or $w$. Then lifting the pair $\{uv,vw\}$ consists of deleting $uv$ and $vw$ and adding a new edge $uw$. In $L(G)$, this corresponds to deleting two vertices and adding a new one in their place, where the new vertex is incident only to $[\delta(u)\cup \delta(w)]- \delta^*(v)$. 
	
	This is equivalent to contracting the edge of $L(G)$ between the vertices $uv$ and $vw$, and then deleting the resulting edges between the new vertex $uw$ and all vertices representing the members of $\delta^*(v)$. Thus for any sequence of operations which transforms $G$ to $H$, the line graph of each intermediate graph is a minor of $L(G)$. Hence $L(H)$ is a minor of $L(G)$, as desired. \hfill $\Box$ \\
	
	We see in the proof that liftings in $G$ correspond to a sort of restricted contraction in $L(G)$. This hints to us that the converse of Proposition \ref{minor} is not true. For a concrete counterexample, consider $G = K_{3,3}$ and note that $L(G)$ contains $K_5$ as a minor. A graph $H$ having $L(H) \cong K_5$ is either a star (possibly having multiple edges) or a triangle with multiple edges. In either case, $H$ contains a vertex of degree greater than three and hence cannot be immersed in $K_{3,3}$.

	It's clear that Proposition \ref{minor} yields two corollaries to the main theorems of our paper, which we will discuss in turn.
	
	\begin{coro} \label{first}
		If $r\geq 4$ and $G$ is a 4-edge-connected graph with $|G| \geq f_{\ref{main}}(r)$, then $L(G)$ contains $L(C_{2,r})$ as a minor.
	\end{coro}
	
	We would like to point out the relationship between this result and that of Oporowski, Oxley, and Thomas \cite{oot}, in which the authors provide a list of four unavoidable minors of 4-connected graphs. Corollary \ref{first} asserts that if a 4-connected graph is known to be the line graph of a large 4-edge-connected graph, then there is just one unavoidable minor. It is worth noting that $L(C_{2,r})$ contains two of the unavoidable graphs of the result in \cite{oot}.

	\begin{coro} \label{final}
		If $t$ and $r$ are positive integers and $G$ is a $(6t-4)$-edge-connected graph with $|G| \geq f_{\ref{ctr}}(t,r)$, then $L(G)$ contains $L(C_{t,r})$ as a minor. 
	\end{coro}
	
	Corollary \ref{final} states that a graph $G$ which is known to be the line graph of a sufficiently large $(6t-4)$-edge-connected graph contains $L(C_{t,r})$ as a minor. Such a graph is necessarily $(6t-4)$-connected. However, it is not the case that every $(6t-4)$-connected line graph is the line graph of a $(6t-4)$-edge-connected graph. It would be desirable if some strengthening of Corollary \ref{final} would apply to all highly-connected line graphs, rather than just those we know to be line graphs of highly edge-connected graphs. The following theorem, the last in this paper, provides such a strengthening. \\
	
	\noindent\textbf{Theorem \ref{linegraph}}
	\textit{For all positive integers $t$ and $r$, there exists a number $f_{\ref{linegraph}}(t,r)$ such that if $G$ is a $(30t-15)$-connected line graph and $|G| \geq f_{\ref{linegraph}}(t,r)$, then $G$ contains $L(C_{t,r})$ as a minor.}\\

	\noindent\textit{Proof:} Suppose that $G = L(H)$ for some graph $H = (V,E)$ and that $|G| \geq 
	f_{\ref{linegraph}}(t,r) := (tr-1)f_{\ref{rooted}}(t,r)$. For ease of notation, let $c(t) = 15t-6$. Note that $30t-15 = 2c(t) - 3$. Let $S = \{v\in V(H): d_H(v) \geq c(t)\}$. 
	
	We first claim that $S$ is $c(t)$-edge-connected in $H$. Otherwise, there exists a minimal edge-cut of the form $\delta(X)$ such that $|\delta(X)| < c(t)$ and $H \bs \delta(X)$ consists of two connected components, each containing members of $S$. Note that because each vertex in $S$ is incident with at least $c(t)$ edges, there must be at least one edge of $H$ in $H[X]$ and $H[V-X]$. This implies the existence of a vertex-cut of size less than $c(t)$ in $G$, a contradiction. 
	
	Next, we claim that $S$ is a vertex cover of $H$; that is, that $H-S$ is edgeless. Suppose otherwise. Then there exists an edge $xy$ of $H$ where $x,y\notin S$. Consider the edge-cut $\delta(\{x,y\})$. Let $p$ denote the number of edges in parallel between $x$ and $y$. Since $x,y\notin S$, we have $d_H(x), d_H(y) \leq c(t) - 1$. This implies that each of $x$ and $y$ contribute at most $c(t) - 1 - p$ edges to $\delta(\{x,y\})$. Hence $|\delta(\{x,y\})| \leq 2c(t) -2 -2p$. Note that $||H|| = |G| > 2c(t) - 3$, and hence the number of edges in $H[V-\{x,y\}]$ is at least $2c(t)-2 - (2c(t)-2 -2p) - p = p$, which is at least 1 by assumption. This implies that $\delta(\{x,y\})$ corresponds to a vertex-cut of $G$ of size at most $2c(t)-4$, a contradiction. 
	
	If $H$ has a vertex of degree $tr$, then this corresponds in $G$ to a $K_{tr}$ subgraph, which in turn implies the presence of a $L(C_{t,r})$ minor in $G$. Hence we may assume that $\Delta(H) < tr$. 
	
	Finally, we find that $$(tr-1)|S| \geq \sum_{v\in S} d_H(v) \geq ||H|| = |G| \geq (tr-1)f_{\ref{rooted}}(t,r)$$ so that $|S| \geq f_{\ref{rooted}}(t,r)$. Theorem \ref{rooted} then implies that $H$ contains $C_{t,r}$ as an immersion, and Corollary \ref{final} implies that $G$ contains $L(C_{t,r})$ as a minor. \hfill $\Box$ \\
	
	We end with a few final comments on this result as it pertains to Theorem \ref{bigktk}. If one's only desire is to force the presence of a large $2t$-connected minor, then Theorem \ref{bigktk} allows them to do so in a graph with connectivity approximately $31t$. Note that $L(C_{t,r})$ is $2t$-connected, so that if one wishes to force the presence of a large $2t$-connected minor in a line graph, then Theorem \ref{linegraph} allows one to do so in a graph of connectivity approximately $30t$, a slight improvement. Additionally, if one cares less about the connectivity function and more about finding as much structure as possible in their graph, we note that $L(C_{t+1,r})$ contains $K_{t,r}$ as a minor, and hence is a ``better" unavoidable minor for the class of line graphs.

\end{document}